\documentclass{amsart}
\usepackage{amscd,amssymb}

\newcommand{\HDD[2]}{{#1}\mbox{-}\overline{\mathrm{DD}\hskip-1pt}\hskip1pt{}^{\{#2\}}}
\newcommand{\DDP[2]}{#1\mbox{-}\mathrm{DD}^{#2}} 
\newcommand{\MAP[2]}{#1\mbox{-}\mathrm{MAP}^{#2}} 
\newcommand{\HDDP[2]}{{#1}\mbox{-}\overline{\mathrm{DD}\hskip-1pt}\hskip1pt{}^{#2}}
\newcommand{\LC[1]}{\mathrm{LC}^{#1}}
\newcommand{\ANE[1]}{\mathrm{ANE}[#1]}
\newcommand{\AP[1]}{\mathrm{AP}[#1]}

\newcommand{\lCH}{\mathcal{LR}}
\newcommand{\HDDmore}{\underset{{\mbox{\tiny DD}}}{\Rightarrow}}
\newcommand{\Tor}{\mathrm{Tor}}
\newcommand{\pTor}{p\mbox{-}\mathrm{Tor}}

\newcommand{\IR}{\mathbb R}
\newcommand{\IN}{\mathbb N}
\newcommand{\IC}{\mathbb C}
\newcommand{\IZ}{\mathbb Z}
\newcommand{\IQ}{\mathbb Q}
\newcommand{\II}{\mathbb I}
\newcommand{\IF}{\mathbb F}
\newcommand{\Rp}{R_p}

\newcommand{\ANR}{\mathrm{ANR}}
\newcommand{\AR}{\mathrm{AR}}
\newcommand{\Br}{\mathsf{Br}}
\newcommand{\trt}{\mathsf{trt}}

\newcommand{\Ddim}{\dim_\triangle}
\newcommand{\DDim}{\dim_\triangle}

\newcommand{\Top}{\mathit{Top}}

\newcommand{\w}{\omega}
\newcommand{\U}{\mathcal U}
\newcommand{\V}{\mathcal V}

\newcommand{\C}{\mathcal C}
\newcommand{\F}{\mathcal F}
\newcommand{\Z}{\mathcal Z}
\newcommand{\DZ}{\overline{\mathcal Z}}

\newcommand{\simU}{\underset{\U}{\sim}}
\newcommand{\pr}{\mathrm{pr}}

\newcommand{\e}{\varepsilon}
\newcommand{\diam}{\mathrm{diam}}

\newcommand{\edim}{\mbox{\rm e-dim}}

\newtheorem{theorem}{Theorem}
\newtheorem{proposition}{Proposition}
\newtheorem{corollary}{Corollary}
\newtheorem{problem}{Problem}

\theoremstyle{definition}
\newtheorem{definition}{Definition}
\newtheorem{remark}{Remark}
\newtheorem{example}{Example}
\newtheorem{question}{Question}

\title{General Position Properties in Fiberwise Geometric Topology}

\author{Taras Banakh}
\address{Department of Mechanics and Mathematics, Ivan Franko 
National University of Lviv (Ukraine) and \newline Instytut Matematyki,
Uniwersytet Humanistyczno-Przyrodniczy Jana Kochanowskiego w Kielcach (Poland)}
\email{tbanakh@yahoo.com}
\thanks{}

\author{Vesko  Valov}
\address{Department of Computer Science and Mathematics, Nipissing University,
100 College Drive, P.O. Box 5002, North Bay, ON, P1B 8L7, Canada}
\email{veskov@nipissingu.ca}
\thanks{The second author was partially supported by NSERC Grant 261914-03.}



\begin{document}
\maketitle

\section*{Introduction}

The classical Lefschetz-Menger-N\"obeling-Pontrjagin-Tolstova
Embedding Theorem asserts that each $n$-dimensional compact metric
space embeds into the \mbox{$(2n+1)$}-dimensional Euclidean space
$\IR^{2n+1}$. A parametric version of this result was proved by
B.Pasynkov \cite{Pas1}: any $n$-dimensional map $p:K\to M$ between
metrizable compacta with $\dim M=m$ embeds into the projection
$\pr_M:M\times\IR^{2n+1+m}\to M$ in the sense that there is an
embedding $e:K\to M\times \IR^{2n+1+m}$ with $\pr_M\circ e=p$. We
recall that  a map $p:X\to Y$ is {\em $n$-dimensional} if $\dim
p^{-1}(y)\le n$ for all $y\in Y$ (\mbox{0-dimensional} maps are
called sometimes {\em light} maps).

The key ingredient of Pasynkov's proof consists in constructing a map
$f:K\to\IR^{2n+1+m}$ that is injective on each fiber $p^{-1}(y)$,
$y\in M$, of $p$. Actually, Pasynkov proved that the set of all such maps
is dense and $G_\delta$ in the function
space $C(K,\IR^{2n+1+m})$.

In this paper we introduce the $m$-DD$^n$-property and show that the
Pasynkov result remains true if $\IR^{2n+1+m}$ is replaced by any
completely metrizable $\LC[m+n]$-space $X$ possessing this property.
Let us note that the $m$-DD$^n$-property is a parametric version of
the classical disjoint $n$-disks property DD$^n$P which plays a
crucial role in characterizing finite and infinite-dimensional
manifolds, see \cite{West}. We shall also give convenient
``arithmetic'' tools for establishing the $m$-DD$^n$-property of
products and obtain on this base simple proofs of some classical and
recent results on (fiber) embeddings. In particular, the Pasynkov
theorem mentioned above, as well as the results of P. Bowers
\cite{Bowers} and Y. Sternfeld \cite{Stern} on embedding into
product of dendrites follow from our general approach. Moreover, the
arithmetics of the $\HDD[m]{n,k}$-properties established in our
paper generalizes some results of W.~Mitchell \cite{mitchell},
R.~Daverman \cite{Dav0}  and D.~Halverson \cite{hal}, \cite{DH}.

\section{Survey of Principal Results}

Throughout the paper $m,n,k$ will stand for non-negative integers or
$\infty$. We extend the arithmetic operations from
$\w=\{0,1,2,\dots\}$ onto $\w\cup\{\infty\}$ letting
$\infty=\infty+\infty=\infty+n=n+\infty=\infty-n$ for any $n\in\w$.
$\mathbb I$ denotes the unit interval $[0,1]$ and  $\IQ$ the set of
rational numbers on the real line $\IR$.  By a simplicial complex we
shall always mean the geometric realization of an abstract
simplicial complex equipped with the $CW$-topology. All topological
spaces are assumed to be {\em Tychonoff} and  all maps continuous.

By an $\ANR$-space we mean a metrizable space $X$ which is a neighborhood retract
of every metrizable space $M$ containing $X$ as a closed subspace.
It is well-known (see \cite{Bors} or \cite{Hu}) that a metrizable
space $X$ is an $\ANR$ if and only if it is an $\ANE[\infty]$ for
the class of metrizable spaces. We recall that a space $X$ is called
an $\ANE[n]$ for a class $\C$ of spaces if every map $f:A\to X$
defined on a closed subset of a space $C\in\C$ with $\dim C\le n$
can be extended to a continuous map $\bar f:U\to X$ defined on some
neighborhood $U$ of $A$ in $X$.

Following \cite{EW}, we define a subset $A$ of a space $X$ to be
{\em relative} $\LC[n]$ in $X$ if given $x\in X$, $k<n+2$, and a
neighborhood $U$ of $x$ there is a neighborhood $V\subset U$ of $x$
such that each map $f:\partial \mathbb I^k\to A\cap V$ extends to a
map $\bar f:\mathbb I^k\to U\cap A$. A space $X$ is an  {\em
$\LC[n]$-space} if it is relative $\LC[n]$ in $X$. According to
\cite[V.2.1]{Hu}, a metrizable space $X$ is $\LC[n]$ for a finite
number $n$ if and only if $X$ is $\ANE[n+1]$ for the class of
metrizable spaces.

\subsection{$\DDP[m]{n}$-property and fiber embeddings}
We recall that a space $X$ has DD$^n$P, the {\em disjoint
$n$-disks property}, if any two maps $f,g:\mathbb I^n\to X$ from the
$n$-dimensional cube $\mathbb I^n=[0,1]^n$ can be approximated by maps
with disjoint images. A parametric version of this property says
that the same can be done for a continuous family $f_z,g_z:\mathbb I^n\to
X$ of maps parameterized by points $z$ of some space $M$.
More precisely, given a compact space $M$, we shall say that a
space $X$ has the {\em $M$-parametric disjoint $n$-disk property}
(briefly, the {\em $M\mbox{\rm -DD}^n$-property}) if any two maps
$f,g:M\times \mathbb I^n\to X$ can be uniformly approximated by maps
$f',g':M\times \mathbb I^n\to X$ such that for any $z\in M$ the images
$f'(\{z\}\times \mathbb I^n)$ and $g'(\{z\}\times \mathbb I^n)$ are disjoint.

We are mostly interested in the particular case of this property
with $M=\mathbb I^m$ being the $m$-dimensional cube. In this case we
write $\DDP[m]{n}$ instead of $\DDP[\mathbb I^m]{n}$. In the
extremal cases when $m$ or $n$ is zero, the $m$-DD$^n$-property
turns out to be very familiar. Namely, the 0-DD$^n$-property is
nothing else but the classical disjoint $n$-disks property, while
the $m$-DD$^0$-property is well-known to specialists in fixed point
and coincidence theories: a space $X$ has the $m$-DD$^0$-property
iff any two maps $f,g:\mathbb I^m\to X$ can be approximated by maps
with disjoint graphs!

It is well known (see \cite{To80} or \cite{Ed}) that all one-to-one
maps from a metrizable $n$-dimensional compactum $K$ into a
completely metrizable $\LC[n-1]$-space $X$ possessing the
DD$^n$P-property form a dense $G_\delta$-set in the function space
$C(K,X)$. Our first principal result is just a parametric version of
this embedding theorem.

\begin{theorem}\label{main1}
A completely metrizable $\LC[m+n]$-space $X$ has the
$\DDP[m]{n}$-property if and only if for every perfect map $p:K\to
M$ between finite-dimensional metrizable spaces with $\dim M\le m$
and $\dim(p)\le n$ the function space $C(K,X)$ contains a dense
$G_{\delta}$-set of maps $f:K\to X$ that are injective on each
fiber $p^{-1}(z)$, $z\in M$.
\end{theorem}

The function space $C(K,X)$ appearing in this theorem is endowed
with the source limitation topology whose neighborhood base at a
given function $f\in C(K,X)$ consists of the sets
$$B_\rho(f,\e)=\{g\in C(K,X):\rho(g,f)<\e\},$$ where $\rho$ runs
over all continuous pseudometrics on $X$ and $\e:K\to(0,\infty)$
runs over continuous positive functions on $K$. Here, the symbol
$\rho(f,g)<\e$ means that $\rho(f(x),g(x))<\e(x)$ for all $x\in K$.
To the best of our knowledge, the notion of source limitation
topology was introduced in the literature (see for example,
\cite{Kr}, \cite{munkers}, \cite{mc}) only for metrizable spaces
$X$. In such a case, for a fixed compatible metric $\rho$ on $X$,
the sets $B_\rho(f,\e)$, $\e\in C(K,(0,\infty)$ and $f\in C(K,X)$,
form a base for a topology $\mathcal T_{\rho}$ on $C(K,X)$. If $K$
is paracompact, then the topology $\mathcal T_{\rho}$ does not
depend on the metric $\rho$ \cite{Kr}. Moreover, $\mathcal T_{\rho}$
has the Baire property provided $K$ is paracompact and $X$ is
completely metrizable \cite{munkers}. According to
Lemma~\ref{func-sp} below, $\mathcal T_{\rho}$ coincides with the
topology obtained from our definition provided $K$ is paracompact
and $X$ metrizable. Therefore, the source limitation topology on
$C(K,X)$ also has the Baire property when $K$ is paracompact and $X$
is completely metrizable. In the sequel, we will use our more
general definition (in terms of pseudometrics) and, unless stated
otherwise, all function spaces will be considered with this
topology.

In fact, finite-dimensionality of the spaces $K,M$ in
Theorem~\ref{main1} can be replaced by the $C$-space property. We
recall that a topological space $X$ is defined to be a {\em
$C$-space} if for any sequence $\{\V_n\colon n\in\omega\}$ of open
covers of $X$ there exists a sequence $\{\U_n\colon n\in\omega\}$ of
disjoint families of open sets in $X$ such that each $\U_n$ refines
$\V_n$ and $\bigcup\{\U_n\colon n\in\omega\}$ is a cover of $X$. It
is known that every finite-dimensional paracompact space (as well as
every hereditarily paracompact countable-dimensional space) is a
$C$-space and normal $C$-spaces are weakly infinite-dimensional, see
\cite[\S6.3]{End}.

\begin{theorem}\label{main1a}
A completely metrizable locally contractible space $X$ has the
$\DDP[m]{n}$-property if and only if for every perfect map $p:K\to
M$ between metrizable $C$-spaces with $\dim M\le m$
and $\dim(p)\le n$ the function space $C(K,X)$ contains a dense
$G_{\delta}$-set of maps $f:K\to X$ that are injective on each
fiber $p^{-1}(z)$, $z\in M$.
\end{theorem}

\subsection{$\Delta$-dimension of maps} There is a natural
temptation to remove the dimensional restrictions on the
spaces $K,M$ from Theorems~\ref{main1} and \ref{main1a}. This indeed can be done if
we replace the usual dimension $\dim(p)$ of the map $p$ with the
so-called $\Delta$-dimension $\Ddim(p)$ (coinciding with $\dim(p)$
for perfect maps $p$ between finite-dimensional metrizable
spaces.)

By definition, the {\em $\Delta$-dimension} $\Ddim(p)$ of a map
$p:X\to Y$ between Tychonoff spaces is equal to the smallest
cardinal number $\tau$ for which there is a map $g:X\to \mathbb I^\tau$
such that the diagonal product $f\Delta g:X\to Y\times \mathbb I^\tau$ is
a light map. The $\Delta$-dimension $\Ddim(p)$ is a well-defined
cardinal function non-exceeding the weight $w(X)$ of $X$ (because
we always can take $g$ to be an embedding in the Tychonoff  cube
$\mathbb I^{w(X)}$).

The following important result describing the interplay between the
dimension and $\Delta$-dimension of perfect maps is actually a
reformulation of results due to B.~Pasynkov \cite{Pas2}, M.~Tuncali,
V.~Valov \cite{TV}, and M.~Levin \cite{Lev}, see
Section~\ref{Deltadim}.

\begin{proposition}\label{Ddim1} Let $f:X\to Y$ be a perfect map between
paracompact spaces. Then
\begin{enumerate}
\item $\dim(f)\le \Ddim(f)$;
\item $\Ddim(f)=0$ if and only if $f$ is a light map;
\item $\DDim(f)\le\w$ if $X$ is submetrizable;
\item $\Ddim(f)=\dim(f)$ if $X$ is submetrizable and $Y$ is a $C$-space;
\item $\DDim(f)\le\dim (f)+1$ if the spaces $X,Y$ are compact and metrizable.
\end{enumerate}
\end{proposition}

We recall that a topological space $X$ is {\em submetrizable} if it
admits a continuous metric (equivalently, admits a bijective
continuous map onto a metrizable space). The following theorem is a
version of Theorem~\ref{main1} with $\dim(p)$ replaced by
$\Ddim(p)$.

\begin{theorem}\label{main2}
A completely metrizable $\ANR$-space $X$ has the
$\DDP[m]{n}$-property if and only if for every perfect map $p:K\to
M$ between submetrizable paracompact spaces with $\dim M\le m$ and
$\Ddim(p)\le n$ the function space $C(K,X)$ contains a dense
$G_{\delta}$-set of maps $f:K\to X$ that are injective on each
fiber $p^{-1}(z)$, $z\in M$.
\end{theorem}

\subsection{The $\HDDP[m]{n}$-property and a general fiber embedding theorem}
In fact, it is more convenient to work not with the
$\DDP[m]{n}$-property, but with its homotopical version defined as
follows:

\begin{definition} A space $X$ has the
{\em $\HDDP[m]{n}$-property} if for any open cover $\U$ of $X$ and
any two maps $f,g:\mathbb I^m\times \mathbb I^n\to X$ there are maps
$f',g':\mathbb I^m\times \mathbb I^n\to X$ such that
\begin{itemize}
\item $f'$ is $\U$-homotopic to $f$;
\item $g'$ is $\U$-homotopic to $g$;
\item $f'(\{z\}\times \mathbb I^n)\cap g'(\{z\}\times \mathbb
I^n)=\varnothing$ for all $z\in \mathbb I^m$.
\end{itemize}
\end{definition}

We recall that two maps $f,g:K\to X$ are said to be {\em
$\U$-homotopic} (briefly, $f\simU g$), where $\U$ is a cover of $X$,
if there is a homotopy $h:K\times[0,1]\to X$ such that for every
$x\in K$ we have $h(x,0)=f(x)$, $h(x,1)=g(x)$ and
$h(\{x\}\times[0,1])$ is contained in some $U\in\U$. It is clear
that any $\U$-homotopic maps $f,g:K\to X$ are {\em $\U$-near} (i.e.,
for each point $z\in K$ the set $\{f(z),g(z)\}$ lies in some
$U\in\U$).

The notion of a $\U$-homotopy has a pseudometric counterpart. Given
a continuous pseudometric $\rho$ on $X$ and a continuous map
$\e:K\to(0,\infty)$ we shall say that two maps $f,g:K\to X$ are
{\em $\e$-homotopic} if there is a homotopy $h:K\times[0,1]\to X$
such that $h(z,0)=f(z)$, $h(z,1)=g(z)$ and
$\diam_\rho(h(\{z\}\times[0,1])<\e(z)$ for all $z\in K$. In this
case $h$ is called an {\em $\e$-homotopy}.

The relation between the $\DDP[m]{n}$-property and its homotopical
version is described by next proposition.

\begin{proposition}\label{lc-HDD} Each space $X$ with the
$\HDDP[m]{n}$-property has the $\DDP[m]{n}$-property. Conversely,
each $\LC[n+m]$-space $X$ possessing the $\DDP[m]{n}$-property has
the $\HDDP[m]{n}$-property.
\end{proposition}

Proposition \ref{lc-HDD} follows from the well-known property of
$\LC[n]$-spaces which asserts that for any open
cover $\U$ of an $\LC[n]$-space $X$ with $n<\infty$ there is
another open cover $\V$ of $X$ such that two maps $f,g:\mathbb I^n\to X$
are $\U$-homotopic provided they are $\V$-near, see Lemma~\ref{lcn1}.

Thus, in the realm of $\LC[m+n]$-spaces both the $\HDDP[m]{n}$-property and the
$\DDP[m]{n}$-property are equivalent.
 The advantage of the $\HDDP[m]{n}$-property is that it works for spaces
 without a nice local structure, while the $m$-DD$^n$-property is
 applicable only for $\LC[k]$-spaces with sufficiently large $k$.
In particular,  using the $\HDDP[m]{n}$-property, we can establish
the following general result
implying Theorems~\ref{main1}, \ref{main1a} and \ref{main2}.

\begin{theorem}\label{main3} Let $p:K\to M$ be a perfect map defined on a
paracompact submetrizable space $K$. If a subspace $X$ of a
completely-metrizable space $Y$ possesses the $\HDDP[m]{n}$-property
for $m=\dim M$ and $n=\Ddim(p)$, then
 $$\mathcal E(p,Y)=\{f\in C(K,Y):\mbox{$p\triangle f:K\to M\times Y$ is an embedding}\}$$
is a $G_\delta$-set in $C(K,Y)$ whose closure $\overline{\mathcal
E(p,Y)}$ contains all simplicially factorizable maps from $K$ to
$X$. More precisely, for any continuous pseudometric $\rho$ on $Y$,
a continuous function $\e:K\to(0,+\infty)$ and a simplicially
factorizable map $f:K\to X$ there is a map $g\in\mathcal E(p,Y)$ and
an $\e$-homotopy $h:K\times[0,1]\to Y$ connecting $f$ and $g$ so
that $h(K\times[0,1))\subset X$.
\end{theorem}

A map $f:K\to X$ is called {\em simplicially factorizable} if there
exist a simplicial complex $L$ and maps $\alpha:K\to L$ and
$\beta:L\to X$ such that $f=\beta\circ\alpha$. It turns out that in
many important cases simplicially factorizable maps form a dense set
in the function space $C(K,X)$. To describe such cases, we need the
notion of a Lefschetz $\ANE[n]$-space that is a parameterized
version of a space satisfying the Lefschetz condition, see
\cite[V.8]{Bors}.

Let $\U$ be a cover of a space $X$ and $K$ be a simplicial complex.
By a {\em partial\/  $\U$-realization} of $K$ in $X$ we understand
any continuous map $f:L\to X$ defined on a geometric subcomplex
$L\subset K$ containing all vertices of $K$ and such that $\diam
f(\sigma\cap L)<\U$ for every simplex $\sigma$ of $K$. If $L=K$,
then the map $f$ is called {\em a full $\U$-realization} of $K$ in
$X$.

A topological space $X$ is defined to be a {\em Lefschetz $\ANE[n]$}
if for every open cover $\U$ of $X$ there is an open cover $\V$ of
$X$ such that each partial $\V$-realization $f:L\to X$ of a
simplicial complex $K$ with $\dim K\le n$ can be extended to a full
$\U$-realization $\tilde f:K\to X$ of $K$.

Lefschetz $\ANE[n]$-spaces are tightly connected with both $\ANR$'s
and $\LC[n]$-spaces and have all basic properties of such spaces.

\begin{proposition}\label{lefschetz} Let $n$ be a non-negative integer or infinity.
\begin{enumerate}
\item A metrizable space $X$ is a Lefschetz $\ANE[n]$ if and only if $X$
is an $\ANE[n]$ for the class of metrizable spaces;
\item If $n$ is finite, then a regular $($paracompact$)$ space $X$ is a
Lefschetz $\ANE[n]$ $($if and$)$ only if $X$ is $\LC[n-1]$;
\item Each convex subset $X$ of a $($locally convex$)$ linear topological space $L$
is a Lefschetz $\ANE[n]$ for any finite $n$ \textup{(}is a Lefschetz
$\ANE[\infty]$\textup{)};
\item There exists a metrizable $\sigma$-compact  linear topological space that fails to be a Lefschetz
$\ANE[\infty]$;
\item A neighborhood retract of a Lefschetz $\ANE[n]$-space is a Lefschetz
$\ANE[n]$-space;
\item A functionally open subspace of a Lefschetz $\ANE[n]$ is a Lefschetz
$\ANE[n]$;
\item A topological space $X$ is a Lefschetz $\ANE[n]$ if $X$ has a uniform open cover by Lefschetz
$\ANE[n]$-spaces;
\item A metric space $(X,\rho)$ is a Lefschetz $\ANE[n]$ if for every $\e>0$ there is $\delta>0$
such that each partial $\mathcal B_\rho(\delta)$-realization $f:L\to
X$ of a simplical complex $K$ with $\dim K\le n$ extends to a full
$\mathcal B_\rho(\e)$-realization $\bar f:K\to X$ of $K$ in $X$;
\item For each continuous pseudometric $\eta$ on a paracompact Lefschetz $\ANE[n]$-space
$X$ there is a continuous pseudometric $\rho\ge\eta$ such that for
every $r\in(0,1/2]$ each partial $\mathcal D_\rho(r/8)$-realization
$f:L\to X$ of a simplical complex $K$ with $\dim K\le n$ extends to
a full $\mathcal D_\rho(r)$-realization $\bar f:K\to X$ of $K$ in
$X$;
\item Each map $f:X\to Y$ from a paracompact Lefschetz $\ANE[n]$-space to a
metrizable space $Y$ factorizes through a metrizable Lefschetz
$\ANE[n]$-space $Z$ in the sense that $f=g\circ h$ for some maps
$h:X\to Z$ and $g:Z\to Y$.
\end{enumerate}
\end{proposition}

Here by $\mathcal D_\rho(\e)$ we denote the cover of a metric space
$(X,\rho)$ by all open sets of diameter $<\e$. With the notion of Lefschetz
$\ANE[n]$-space at hands, we can returns to simplicially
factorizable maps.

\begin{proposition}\label{simpfactorizable} Simplicially
factorizable maps from a paracompact space $K$ into a Tychonoff
space $X$ form a dense set in the function space $C(K,X)$ if one of
the following conditions is satisfied:
\begin{enumerate}
\item $X$ is a Lefschetz $\ANE[k]$ for $k=\dim K$;
\item $K$ is a $C$-space and $X$ is a locally
contractible paracompact space.
\end{enumerate}
\end{proposition}

Observe that Theorem~\ref{main1}, \ref{main1a}, and \ref{main2}
follow immediately from Theorem~\ref{main3} and
Propositions~\ref{simpfactorizable}, \ref{lefschetz} and
\ref{Ddim1}.

Combining  Theorem~\ref{main3} with
Propositions~\ref{simpfactorizable}(1), \ref{lefschetz}(2) and
\ref{Ddim1}(4), we obtain another generalization of
Theorem~\ref{main1}.

\begin{theorem} \label{main4} Let $p:K\to M$ be a perfect map
between  finite-dimensional paracompact spaces with $K$ being
submetrizable. If $X$ is a completely metrizable $\LC[k-1]$-space
possessing the $\HDDP[m]{n}$-property, where $k=\dim K$, $m=\dim M$
and $n=\dim(p)$, then the function space $C(K,X)$ contains a dense
$G_\delta$-set of maps injective on each fiber of $p$.
\end{theorem}

\subsection{Approximating perfect maps by perfect PL-maps}

The proof of Theorem~\ref{main3} heavily exploits the technique of
approximations by PL-maps. By a {\em PL-map} (resp., a {\em
simplicial map}) we understand a map $f:K\to M$ between simplicial
complexes which maps each simplex $\sigma$ of $K$ into (resp., onto)
some simplex $\tau$ of $M$ and $f$ is linear on $\sigma$.

\begin{theorem}\label{PL} If $p:X\to Y$ is a perfect map between
paracompact spaces, then for any open cover $\U$ of $X$ there exists
an open cover $\V$ of $Y$ such that for any $\V$-map $\beta:Y\to M$
into a simplicial complex $M$ there are an $\U$-map $\alpha:X\to K$
into a simplicial complex $K$ and a perfect PL-map $f:K\to M$ with
$f\circ\alpha=\beta\circ p$ and $\Ddim(f)=\dim(f)\le\Ddim(p)$.
\end{theorem}

Since for each open cover $\V$ of a paracompact space $Y$ there is a
$\V$-map $\beta:Y\to M$ into a simplicial complex of dimension $\dim
M\le\dim Y$, Theorem~\ref{PL} implies the following approximation
result.

\begin{corollary}\label{PL:cor} If $p:X\to Y$ is a perfect map between
paracompact spaces, then for any open covers $\U$ and $\V$ of $X$
and $Y$, respectively, there exist a $\U$-map $\alpha:X\to K$ into
a simplicial complex $K$ of dimension $\dim K\le\dim Y+\Ddim(p)$, a
$\V$-map $\beta:Y\to M$ to a simplicial complex $M$ of dimension
$\dim M\le\dim Y$, and a perfect PL-map $f:K\to M$ of dimension
$\dim(f)\le\Ddim(p)$ making the following diagram commutative: $$
\begin{CD}
X@>{\alpha}>>K\cr @V{p}VV @VV{f}V\cr Y@>{\beta}>>M\cr
\end{CD}
$$
\end{corollary}

For light maps $p:X\to Y$ between metrizable compacta this
corollary was proved by A.Dranishnikov and V.Uspenskij in
\cite{DU} and for arbitrary maps between metrizable compacta by
Yu.Turygin \cite{Tu}.

\subsection{$\HDD[m]{n,k}$-properties}
Because of the presence of the $\HDDP[m]{n}$-property in
Theorems~\ref{main3}--\ref{main4}, it is important to have
convenient methods for detecting that property. To establish such methods,
we introduce the following three-parametric version of
$\HDDP[m]{n}$.

\begin{definition} A space $X$ is defined to have the {\em
$\HDD[m]{n,k}$-property} if for any open cover $\U$ of $X$ and two
maps $f:\mathbb I^m\times \mathbb I^n\to X$, $g:\mathbb I^m\times
\mathbb I^k\to X$ there exist maps $f':\mathbb I^m\times \mathbb
I^n\to X$, $g':\mathbb I^m\times \mathbb I^k\to X$ such that
$f'\simU f$, $g'\simU g$, and $f'(\{z\}\times \mathbb I^n)\cap
g'(\{z\}\times \mathbb I^k)=\varnothing$ for all $z\in \mathbb I^m$.
\end{definition}

By $\HDD[m]{n,k}$ we shall denote the class of all spaces with the
$\HDD[m]{n,k}$-property. It is clear that the $\HDDP[m]{n}$-property
coincides with the $\HDD[m]{n,n}$-property. If some of the numbers
$m,n,k$ are infinite, the detection of the $\HDD[m]{n,k}$-property
can be reduced to the detection of $\HDD[m]{n,k}$ with finite
$m,n,k$.

\begin{proposition}\label{infreduction} A Tychonoff  space $X$ has
the $\HDD[m]{n,k}$-property if and only if it has the
$\HDD[a]{b,c}$-property for all $a<m+1$, $b<n+1$, $c<k+1$.
\end{proposition}

The proof of Theorem~\ref{main3} is based on the following
simplicial characterization of the $\HDD[m]{n,k}$-property.

\begin{theorem}\label{simplicial} A submetrizable space $X$ has the
$\HDD[m]{n,k}$-property if and only if for any
\begin{itemize}
\item  simplicial maps $p_N:N\to M$, $p_K:K\to M$ between
finite simplicial complexes with $\dim M\le m$, $\dim(p_N)\le n$,
$\dim(p_K)\le k$,
\item  open cover $\U$ of $X$, and
\item maps $f:N\to X$, $g:K\to X$,
\end{itemize}
there exist maps $f':N\to X$, $g':K\to X$ such that $f'\simU f$,
$g'\simU g$ and, for every $z\in M$ we have
\mbox{$f'(p_N^{-1}(z))\cap g'(p_K^{-1}(z))=\varnothing$}.
\end{theorem}

Using the above simplicial characterization, we can establish the
local nature of the $\HDD[m]{n,k}$-property.

\begin{proposition}\label{local} Let $m,n,k$ be non-negative
integers or $\infty$.
\begin{enumerate}
\item If a space $X$ has the $\HDD[m]{n,k}$-property, then each
open subspace of $X$ also has that property.
\item A paracompact
submetrizable space $X$ has the $\HDD[m]{n,k}$-property if and only
if it admits a cover by open subspaces with that property.
\end{enumerate}
\end{proposition}

The $\HDD[m]{n,k}$-property is also preserved by taking
homotopically $n$-dense subspaces. We define a subset $A$ of a
topological space $X$ to be {\em homotopically $n$-dense} in $X$ if
the following conditions are satisfied:
\begin{itemize}
\item for every map $f:\mathbb I^n\to X$ and an open cover $\U$ of
$X$ there is a map $f':\mathbb I^n\to A$ that is $\U$-homotopic to
$f$;
\item for every open cover $\U$ of $X$ there is an open cover
$\V$ of $X$ such that if two maps $f,g:\mathbb I^{n}\to A$  are
$\V$-homotopic in $X$, then they are $\U$-homotopic in $A$.
\end{itemize}
By Theorem 2.8 of \cite{To78}, each dense relative $\LC[n]$-subset
$X$ of a metrizable space $\tilde X$ is homotopically $n$-dense in
$\tilde X$. The following useful proposition follows immediately
from the definitions and the mentioned theorem of Toru\'nczyk.

\begin{proposition}\label{negligible} A homotopically
$\max\{m+n,m+k\}$-dense subspace $X$ of a topological space $\tilde
X$ has the $\HDD[m]{n,k}$-property if and only if $\tilde X$ has
that property. Consequently, a dense relative
$\LC[m+\max\{n,k\}]$-set $X$ in a space $\tilde X$ has the
$\HDD[m]{n,k}$-property if and only if $\tilde X$ has that property.
\end{proposition}

This fact will be often applied in combination with Proposition 2.8
from \cite{DM} asserting that each metrizable $\LC[n]$-space $X$
embeds into a completely metrizable $\LC[n]$-space $\tilde X$ as a
dense relative $\LC[n]$-set.
 The last assertion enable us  to apply Baire Category
arguments for establishing the $\HDD[m]{n,k}$-properties in
arbitrary (not necessary complete) metric spaces.

Next, we elaborate tools for detecting the
$\HDD[m]{n,k}$-properties. Recall that a space $X$ {\em has no free
arcs} if $X$ contains no open subset, homeomorphic to a non-empty
connected subset of the real line. In particular, a space without
free arcs has no isolated points.

\begin{proposition}\label{lc-00}
\begin{enumerate}
\item A topological space $X$ has the $\HDD[0]{0,0}$-property if
and only if each path-connected component of $X$ is non-degenerate.
\item  An $\LC[0]$-space $X$ has the $\HDD[0]{0,0}$-property if and
only if $X$ has no isolated point.
\item A metrizable $\LC[1]$-space
$X$ has the $\HDD[0]{0,1}$-property iff $X$ has the
$\HDD[1]{0,0}$-property iff $X$ has no free arc.
\item Any
metrizable $\LC[n]$-space $X$ with the $\HDD[0]{0,n}$-property and
$\dim X\le n$ has the $\HDD[0]{0,\infty}$-property.
\item A Polish
$\ANE[\max\{n,k\}+1]$-space $X$ has the $\HDD[0]{n,k}$-property if and
only if there are two disjoint dense $\sigma$-compact subsets $A,B$ of $X$
such that $A$ is relative $\LC[n-1]$ and $B$ is relative $\LC[k-1]$ in $X$.
\end{enumerate}
\end{proposition}

Items 3 and 4 of Proposition~\ref{lc-00} imply that each one dimensional
$\LC[1]$-space without free arcs has the $\HDD[0]{0,\infty}$-property.
In particular, each dendrite with a dense set of end-points has that
property.

The last item of Proposition~\ref{lc-00} is a partial case of a more
general characterization of the $\HDD[m]{n,k}$-property in terms of
mapping absorption properties.

Let $M,X$ be topological spaces. We shall say that a subset
$A\subset M\times X$ has the {\em absorption property for
$n$-dimensional maps in $M$} (briefly, $\MAP[M]{n}$) if for any
$n$-dimensional map $p:K\to M$ with $K$ being a finite-dimensional
compact space, a closed subset $C\subset K$, a map $f:K\to X$, and
an open cover $\U$ of $X$ there is a map $f':K\to X$ such that $f'$
is $\U$-homotopic to $f$, $f'|C=f|C$ and $(p\triangle f')(K\setminus
C)\subset A$. If $M=\mathbb I^m$, then we write $\MAP[m]{n}$ instead
of $\MAP[\mathbb I^m]{n}$.

\begin{theorem}\label{map} Let $m,n,k$ be non-negative integers or
infinity and $d=1+m+\max\{n,k\}$. A \textup{(}Polish
$\ANE[d]$-\textup{)}space $X$ has the $\HDD[m]{n,k}$-property if
\textup{(}and only if\textup{)} for any separable polyhedron $M$
with $\dim M\le m$ there are two disjoint
\textup{(}$\sigma$-compact\textup{)} sets $E,F\subset M\times X$
such that $E$ has $\MAP[M]{n}$ and $F$ has $\MAP[M]{k}$.
\end{theorem}

Let us observe that the existence of such disjoint sets $A,B$ is not
obvious even for a dendrite with a dense set of end-points. Such a
dendrite $D$ has the $\HDD[1]{0,0}$-property and thus the product
$\mathbb I\times D$ contains two disjoint $\sigma$-compact subsets
with $\MAP[1]{0}$.

\subsection{A Selection Theorem for $Z_n$-set-valued functions}
Many results on $\HDD[m]{n,k}$-properties are based on a selection
theorem for $Z_n$-valued functions, discussed in this subsection.

A subset $A$ of a topological space $X$ is called a ({\em
homotopical\/}) $Z_n$-set in $X$ if $A$ is closed in $X$ and for any
an open cover $\U$ of $X$  and a map $f:\II^n\to X$ there is a map
$g:\II^n\to X$ such that $g(\II^n)\cap A=\varnothing$ and $g$ is
 $\U$-near ($\U$-homotopic) to $f$.
Each homotopical $Z_n$-set in a topological space $X$ is a $Z_n$-set in $X$.
The converse is true if $X$ is an $\LC[n]$-space, see Theorem~\ref{Z2+hZn=Zn}.

A set-valued function $\Phi:X\Rightarrow Y$ is defined to be {\em
compactly semi-continuous} if for every compact subset $K\subset Y$
the preimage $\Phi^{-1}(K)=\{x\in X:\Phi(x)\cap K\ne\varnothing\}$
is closed in $X$.

\begin{theorem}\label{selection} Let $\Phi:X\Rightarrow Y$ be a compactly semi-continuous
set-valued function from a paracompact $C$-space $X$ into a
topological space $Y$, assigning to each point $x\in X$ a
homotopical $Z_n$-set $\Phi(x)$, where $n=\dim X\leq\infty$. If $X$
is a retract of an open subset of a locally convex linear
topological space, then for any map $f:X\to Y$ and any continuous
pseudometric $\rho$ on $Y$ there is map $f':X\to Y$ such that
$f'(x)\notin\Phi(x)$ for all $x\in X$ and $f'$ is $1$-homotopic to
$f$ with respect to $\rho$.
\end{theorem}

In particular, the theorem is true for stratifiable $\ANR$'s $X$
(which are neighborhood retracts of stratifiable locally convex
spaces, see \cite{Sipa}).

\subsection{Homotopical $Z_n$-sets and $\HDD[m]{n,k}$-properties}
It turns out that homotopical $Z_n$-sets are tightly connected with
the $\HDD[m]{n,k}$-properties. A point $x$ of a topological space
$X$ is called a ({\em homotopical}\/) {\em $Z_n$-point} if the
singleton $\{x\}$ is a (homotopical) $Z_n$-set in $X$. By $\Z_n(X)$
we shall denote the set of all homotopical $Z_n$-points of a space
$X$.

Let
\begin{itemize}
\item $\Z_n$ be the class of Tychonoff spaces $X$ with $\Z_n(X)=X$;
\item $\DZ_n$ be the class of Tychonoff spaces $X$ with $\overline{\Z_n(X)}=X$;
\item $\Delta\Z_n$ be the class of Tychonoff spaces $X$ whose
diagonal $\Delta_X$ is a homotopical $Z_n$-set in $X^2$;
\end{itemize}

For example, $\IR^{n+1}\in \Z_n\cap \Delta\Z_n$.

Besides the classes of spaces related to $Z$-sets, we also need some
other (more familiar) classes of topological spaces:
\begin{itemize}
\item $\Br$, the class of metrizable separable Baire spaces,
\item $\Pi^0_2$, the class of Polish spaces, and
\item $\LC[n]$, the class of all $\LC[n]$-spaces.
\end{itemize}

\begin{theorem}\label{Zset} Let $m,n,k$ be non-negative integers or infinity.
\begin{enumerate}
\item A space $X$ has the $\HDD[n]{0,0}$-property if and only if
the diagonal of $X^2$ is a homotopical $Z_n$-set in $X^2$. This can
be written as
\smallskip

\hskip-10pt\frame{\phantom{$I^{I^{I^{I^i}}}$}\hskip-15pt
$\Delta \Z_n=\HDD[n]{0,0}$\hskip-15pt\phantom{$I_{I_{I_{I_I}}}$}}
\smallskip

\item An $\LC[0]$-space $X$ has the $\HDD[0]{0,n}$-property
provided the set $\Z_n(X)$ is dense in $X$. This can be written as:
\smallskip

\hskip-10pt\frame{\phantom{$I^{I^{I^{I^i}}}$}\hskip-15pt
$\LC[0]\cap \DZ_n\subset\HDD[0]{0,n}$\hskip-15pt\phantom{$I_{I_{I_{I_I}}}$}}
\medskip

\item If a metrizable separable Baire \textup{(}$\LC[n]$-\textup{)}space $X$ has the
$\HDD[0]{0,n}$-property then the set of $($homotopical$)$
$Z_n$-points is a dense $G_\delta$-set in $X$:
\smallskip

\hskip-10pt\frame{\phantom{$I^{I^{I^{I^i}}}$}\hskip-15pt
$\Br\cap \LC[n]\cap \HDD[0]{0,n}\subset \DZ_n$\hskip-15pt\phantom{$I_{I_{I_{I_I}}}$}}
\medskip

\item If each point of a space $X$ is a homotopical
$Z_{m+k}$-point, then $X$ has the $\HDD[m]{0,k}$-property:
\smallskip

\hskip-10pt\frame{\phantom{$I^{I^{I^{I^I}}}$}\hskip-15pt
$\Z_{m+k} \subset\HDD[m]{0,k}$\hskip-15pt\phantom{$I_{I_{I_{I_I}}}$}}
\smallskip

\item If a topological space $X$ has either the $\HDD[n]{n,0}$- or
the $\HDD[0]{n,n}$-property, then each point of $X$ is a homotopical
$Z_n$-point:
\smallskip

\hskip-10pt \frame{\phantom{$I^{I^{I^{I^I}}}$}\hskip-15pt
$\HDD[0]{n,n}\cup\HDD[n]{n,0}\subset \Z_n$
\hskip-15pt\phantom{$I_{I_{I_{I_I}}}$}}
\smallskip

\item If a Tychonoff space $X$ has the $\HDD[2]{0,0}$-property,
then each point of $X$ is a homotopical $Z_1$-point:
\smallskip

\hskip-10pt\frame{\phantom{$I^{I^{I^{I^I}}}$}\hskip-15pt
$\HDD[2]{0,0}\subset \Z_1$
\hskip-15pt\phantom{$I_{I_{I_{I_I}}}$}}
\end{enumerate}
\end{theorem}

\subsection{Arithmetics of $\HDD[m]{n,k}$-properties}
In this subsection we study the behavior of the
$\HDD[m]{n,k}$-properties under arithmetic operations. The
combination of the results from this subsection and
Propositions~\ref{local}--\ref{lc-00}  provides convenient tools for
detecting the $\HDD[m]{n,k}$-properties of more complex spaces (like
products or manifolds).

For a better visual presentation of our subsequent results, let us
introduce the following operations on subclasses $\mathcal
A,\mathcal B\subset\Top$ of the class $\Top$ of topological spaces:
$$\begin{gathered}
\mathcal A\times\mathcal B=\{A\times B:A\in\mathcal A,\; B\in\mathcal B\},\\
\frac{\mathcal A}{\mathcal B}=\{X\in\Top:\exists B\in\mathcal B\mbox{ with } X\times B\in\mathcal A\},\\
\mathcal A^k=\{A^k:A\in\mathcal A\}\mbox{  and }\sqrt[k]{\mathcal A}=\{A\in\Top:A^k\in\mathcal A\}.
\end{gathered}
$$
A space $X$ will be identified with the one-element class $\{X\}$.
So $X\times\mathcal A$ and $\dfrac{\mathcal A}X$ means
$\{X\}\times\mathcal A$ and $\dfrac{\mathcal A}{\{X\}}$.

We recall that $\HDD[m]{n,k}$ stands for the class of all spaces
possessing the $\HDD[m]{n,k}$-property and $\LC[n]$ is the class of
$\LC[n]$-spaces.

\begin{theorem}[\bf Multiplication Formulas]\label{arithmetic} Let
$X,Y$ be metrizable spaces and $k_1,k_2,k,\;n_1$, $n_2$,
$n,\;m_1,m_2,m$ be non-negative integers or infinity.
\begin{enumerate}
\item {\bf (First Multiplication Formula)}
\newline
If $X$ has the $\HDD[m]{n,k_1}$-property and $Y$ has the
$\HDD[m]{n,k_2}$-property, then the product $X\times Y$ has the
$\HDD[m]{n,k_1+k_2+1}$-property. This can be written as
\smallskip

\hskip-10pt\frame{\phantom{$I^{I^{I^{I^i}}}$}\hskip-15pt
$\HDD[m]{n,k_1} \times \HDD[m]{n,k_2}\subset\HDD[m]{n,k_1+k_2+1}$\hskip-15pt\phantom{$I_{I_{I_{I_I}}}$}}
\medskip

\item {\bf (Second Multiplication Formula)}
\newline
If $X$ has the $\HDD[m]{n_1,k_1}$-property and $Y$ has both the
$\HDD[m]{n,k_2}$- and $\HDD[m]{n_2,k}$-properties, where
$n=n_1+n_2+1$ and $k=k_1+k_2+1$, then the product $X\times Y$ has
the $\HDD[m]{n,k}$-property. This can be written as
\smallskip

\hskip-10pt\frame{\phantom{$I^{I^{I^{I^I}}}$}\hskip-15pt
$\HDD[m]{n_1,k_1} \times \big(\HDD[m]{n,k_2}\cap\HDD[m]{n_2,k}\big)\subset\HDD[m]{n,k}$
\hskip-15pt\phantom{$I_{I_{I_{I_I}}}$}}
\medskip
\item {\bf (Multiplication by a cell)}\newline  If $X$ has the $\HDD[m]{n,k}$-property, then for any $d<m+1$ the
product $\mathbb I^d\times X$ has the $\HDD[(m-d)]{d+n,d+k}$-property. This can be written as
\smallskip

\hskip-10pt\frame{\phantom{$I^{I^{I^{I^I}}}$}\hskip-15pt $\mathbb I^d\times \HDD[m]{n,k}\subset \HDD[(m-d)]{d+n,d+k}$
\hskip-15pt \phantom{$I_{I_{I_{I_I}}}$}}
\end{enumerate}
\end{theorem}

\begin{remark} Let us mention that, since $\mathbb
R\in\HDD[0]{0,0}$, the second multiplication formula implies the
following result of W. Mitchell \cite[Theorem 4.3(3)]{mitchell} (see
also R. Daverman \cite[Proposition 2.8]{Dav0}): If $X$ is a compact
metric $\ANR$-space with $X\in\HDD[0]{p,p+1}$, then $X\times\mathbb
R\in\HDD[0]{p+1,p+1}$. Moreover, Theorem~\ref{arithmetic}(2) yields
$X\times\mathbb R^{m+p+1}\in\HDD[m]{n+p+1,k}$ for any metrizable
space $X\in\HDD[m]{n,k}\cap\HDD[m]{n+p+1,k-1}$. The partial case of
this result when $m=0$ and $p=1$ was provided  by W. Mitchell in
\cite[Theorem 4.3(2)]{mitchell}. Similarly, we can see that
Theorem~\ref{arithmetic}(3) generalizes the following result of D.
Halverson \cite{hal}: If $X$ is a locally compact {\it ANR} with
$X\in\HDD[1]{1,1}$, then $X\times\mathbb R\in\HDD[0]{2,2}$.
\end{remark}

Next, we consider the so-called base enlargement formulas expressing
the $\HDD[m]{n,k}$-property via $\HDD[0]{n',k'}$-properties for
sufficiently large $n',k'$.

\begin{theorem}[\bf Base Enlargement Formulas]\label{basenlarge}
Let $X$ be a metrizable space and $n,k,m$, $m_1$, $m_2$ be
non-negative integers or infinity.
\begin{enumerate}
\item  If $X$ possesses the \mbox{$\HDD[0]{n+m_1,k+m_2}$-},
$\HDD[m_1]{n,k+m-m_1}$-, and\newline $\HDD[m_2]{n+m-m_2,k}$-properties
simultaneously with $m=m_1+m_2+1$, then $X$ has the
$\HDD[m]{n,k}$-property. This can be written as
\smallskip

\hskip-30pt
\frame{\phantom{$I^{I^{I^{I^I}}}$} \hskip-20pt
$\HDD[0]{n+m_1,k+m_2}\cap\HDD[m_1]{n,k+m-m_1}\cap\HDD[m_2]{n+m-m_2,k}\subset\HDD[m]{n,k}$
\hskip-20pt\phantom{$I_{I_{I_{I_I}}}$}}

\item If $X\in\HDD[0]{n,k+m+1}\cap\HDD[m]{n+1,k}$, then $X$ has the
$\HDD[(m+1)]{n,k}$-property. This can be written as
\smallskip

\hskip-10pt\frame{\phantom{$I^{I^{I^{I^I}}}$}\hskip-15pt
$\HDD[0]{n,k+m+1}\cap\HDD[m]{n+1,k}\subset\HDD[(m+1)]{n,k}$
\hskip-15pt\phantom{$I_{I_{I_{I_I}}}$}}

\item $X$ has the $\HDD[m]{n,k}$-property if $X$ has the
$\HDD[0]{n+i,k+j}$-property for all $i,j\in\w$ with $i+j<m+1$. This
can be written as
\smallskip

\hskip-10pt\frame{\phantom{$I^{I^{I^{I^I}}}$}\hskip-15pt
$\bigcap_{i+j<m+1}\HDD[0]{n+i,k+j}\subset\HDD[m]{n,k}$\hskip-15pt\phantom{$I_{I_{I_{I_I}}}$}}
\end{enumerate}
\end{theorem}

The second base enlargement formula implies that if $X$ is a
metrizable space with $X\in\HDD[0]{1,2}$, then $X\in\HDD[1]{1,1}$.
This result was established by D. Halverson \cite{hal} in the
particular case when $X$ is a separable locally compact {\it ANR}.

\subsection{$\HDD[m]{n,k}$-properties of products} In this
subsection we apply the arithmetic formulas from previous subsection
to establish the $\HDD[m]{n,k}$-properties of products.

\begin{theorem}\label{power} Let
$m,n,k,d,l$ be non-negative integers and $L,D$ be metrizable spaces
such that $L$ has the $\HDD[0]{0,0}$-property and $D$ has the
$\HDD[0]{0,d+l}$-property. If $m+n+k<2d+l$, then the product
$D^d\times L^l$ has the $\HDD[m]{n,k}$-property. This can be written
as
\smallskip

\hskip10pt\frame{\phantom{$I^{I^{I^{I^I}}}$}\hskip-15pt
$\big(\HDD[0]{0,0}\big)^l\times
\big(\HDD[0]{0,d+l}\big)^d\subset\bigcap
\limits_{m+n+k<2d+l}\HDD[m]{n,k}$\hskip-15pt\phantom{$I_{I_{I_{I_{I_I}}}}$}}
\end{theorem}

Combining Theorem~\ref{power} with Theorem~\ref{main3},
Proposition~\ref{lc-HDD} and Proposition~\ref{lc-00}  we obtain

\begin{theorem}\label{power2} Let $l,d$ be non-negative integers or
infinity and $L$, $D$ be completely metrizable locally
path-connected spaces such that $L$ has no isolated points and $D$
is $1$-dimensional without free arcs. Then the product $D^d\times
L^l$ has the $\HDDP[m]{n}$-property for all $m,n\in\w$ with
$m+2n<l+2d$. Consequently, if $p:K\to M$ is a perfect map between
paracompact submetrizable spaces with $\dim M+2\Ddim(p)<l+2d$, then
any simplicially factorizable map $f:K\to D^d\times L^l$ can be
approximated by maps injective on each fiber
 of $p$.
\end{theorem}

The case $m=0$ from  Theorem~\ref{power2} yields

\begin{corollary}\label{cor:power} Let $L,D$ be completely
metrizable $ANR$'s such that $L$ has no isolated points and $D$ is
$1$-dimensional without free arcs. Then the product $D^d\times L^l$
has the $\mathrm{DD}^n{\mathrm P}$ for all $\displaystyle
n<d+\frac{l}2$. Consequently, for any compact space $X$ of dimension
$\displaystyle\dim X<d+\frac{l}{2}$ the set of all embeddings is
dense $G_\delta$ in the function space\/ $C(X,D^d\times L^l)$.
\end{corollary}

\begin{remark} Corollary~\ref{cor:power} generalizes many (if not all)
results on embeddings into products. Indeed, letting $L=\IR$ to be
the real line and $D$ to be a dendrite with a dense set of end-points
we obtain the following well known results:
\begin{itemize}
\item the case $d=0$ and $l=2n+1$ is the
Lefschetz-Menger-N\"obeling-Pontrjagin embedding theorem that
$\IR^{2n+1}$ has $\mathrm{DD}^n\mathrm{P}$;
\item the case $d=n+1$ and $l=0$ is the embedding theorem of P.~Bowers \cite{Bowers}
that $D^{n+1}$ has $\mathrm{DD}^n\mathrm{P}$;
\item the case $d=n$ and $l=1$ is the embedding theorem of Y.~Sternfeld \cite{Stern}
that $D^n\times\mathbb I$ has $\mathrm{DD}^n\mathrm{P}$;
\end{itemize}
Also, for $d=0$ and $m=0$ Theorem~\ref{power2} is close to the
embedding theorem from T.~Banakh, Kh.~Trushchak \cite{BTr} while for
$l=0$ and $m=0$ it is close to that one of T.~Banakh, R.~Cauty,
Kh.~Trushchak, L.~Zdomskyy \cite{BCTZ}.
\end{remark}

\begin{remark} Letting $d=0$ and $L=\IR$ in Theorem~\ref{power2},
we obtain the Pasynkov's result \cite{Pas1} asserting that for a map
$p:X\to Y$ between finite-dimensional metrizable compacta the
function space $C(X,\IR^{\dim Y+2\dim(p)+1})$ contains a dense
$G_\delta$-set of maps that are injective on each fiber of the map
$p$.

Therefore, Theorem~\ref{power2} can be considered as a
generalization of \cite{Pas1}. However, Theorem~\ref{power2} does
not cover another generalization of the Pasynkov's result due to
H.Toru\'nczyk \cite{ht}: {\em If $p:X\to Y$ is a map between
compacta, then the space $C(X,\IR^{\dim X+\dim(p)+1})$ contains a
dense $G_\delta$-set of maps that are injective on each fiber of the
map $p$.}

Taking into account that the Euclidean space $\IR^d$ has the
$\HDD[m]{n,k}$-properties for all $m,n,k$ with $m+n+k<d$, we may ask
whether the this theorem of H.Toru\'nczyk is true in the following
more general form.
\end{remark}

\begin{problem} Let $p:K\to M$ be a map between finite-dimensional compact metric spaces
and $X$ be a Polish $AR$-space possessing
 the $\HDD[m]{n,k}$-property for all $m,n,k$ with $m+n+k\le
\dim K+\dim(p)$. Does $p$ embed into the projection $\pr:M\times
X\to M$\/?
\end{problem}

Let us also note that the mentioned above result of H.Toru\'nczyk would follow from our
Theorem~\ref{main1} if the
 following problem had an affirmative answer.

\begin{problem} Let $f:X\to Y$ be a $k$-dimensional map between
finite-dimensional metrizable compacta. Is it true that there is a
map $g:Y\to Z$ to a compact space $Z$ with $\dim Z\le\dim X-k$
such that the map $g\circ f$ is still  $k$-dimensional\/?
\end{problem}

\subsection{A short survey on homological $Z_n$-sets} The most
exciting results on $\HDD[m]{n,k}$-properties (like multiplication
and $k$-root formulas) are obtained by using homological $Z_n$-sets.
In this subsection we survey some basic facts about such sets, and
refer the interested reader to \cite{BKV} where all these results
are established. We use the singular homology with coefficients in
an Abelian group $G$. If $G=\IZ$, we write $H_k(X)$ instead of
$H_k(X;\IZ)$. By $\widetilde H_*(X;G)$ we denote the singular
homology of $X$ reduced in dimension zero.
\smallskip

It can be shown that a closed subset $A$ of a topological space $X$
is a homotopical $Z_n$-set in $X$ if and only if for every open set
$U\subset X$ the inclusion $U\setminus A\to U$ is weak homotopy
equivalence, which means that the relative homotopy groups
$\pi_k(U,U\setminus A)$ vanish for all $k<n+1$. Replacing the
relative homotopy groups by relative homology groups, we obtain the
notion of a homological $Z_n$-set.

A closed subset $A$ of a space $X$ is defined to be
\begin{itemize}
\item a {\em $G$-homological $Z_n$-set} in $X$ for a coefficient
group $G$ if $H_k(U,U\setminus A;G)=0$ for all open sets $U\subset
X$ and all $k<n+1$;
\item an {\em $\exists G$-homological $Z_n$-set}
if $A$ is a $G$-homological $Z_n$-set in $X$ for some coefficient
group $G$;
\item a {\em homological $Z_n$-set} if $A$ is a
$\IZ$-homological $Z_n$-set in $X$ (equivalently, if $A$ is a
$G$-homological $Z_n$-set for every coefficient group $G$).
\end{itemize}
In \cite{DW} homological $Z_\infty$-sets are referred to as closed
sets of infinite codimension. On the other hand, the term
``homological $Z_n$-set'' has been used in \cite{RS}, \cite{BKV},
\cite{BC}, and \cite{BR}.

The following theorem whose proof can be found in \cite[Theorems
3.2-3.3]{BKV} describes the interplay between various sorts of
$Z_n$-sets.

\begin{theorem}\label{Z2+hZn=Zn} Let $X$ be a topological space.
\begin{enumerate}
\item Each homotopical $Z_n$-set in $X$ is both a $Z_n$-set and a homological $Z_n$-set;
\item Each $Z_n$-set in an $\LC[n]$-space is a homotopical $Z_n$-set;
\item A set is a homotopical $Z_0$-set in $X$ iff it is a $\exists G$-homological $Z_0$-set;
\item Each $\exists G$-homological $Z_1$-set in $X$ is a $Z_1$-set;
\item If $X$ is an $\LC[1]$-space, then a homotopical $Z_2$-set in
$X$ is a homotopical $Z_n$-set if and only if it is a homological
$Z_n$-set.
\end{enumerate}
\end{theorem}

The last item of this theorem has fundamental importance since it
allows application of powerful tools of Algebraic Topology for
studying homotopical $Z_n$-sets and related
$\HDD[m]{n,k}$-properties.

The study of $G$-homological $Z_n$-sets for an arbitrary group $G$
can be reduced to considering Bockstein groups. Under the notations
below,
\begin{itemize}
\item $\IQ$ - the group of rational numbers;
\item $\IZ_p=\IZ/\IZ_p$  - the cyclic group of a prime order $p$;
\item $\IQ_p=\{z\in\IC:\exists k\in\IN \; z^{p^k}=1\}$ - the quasicyclic $p$-group;
\item $\Rp=\{\frac{m}{n}:m\in\IZ$ and $n\in\IN$ is not divisible by $p\}$,
\end{itemize}
the Bockstein family $\sigma(G)$ of a group $G$ is a subfamily of
$\{\IQ,\IZ_p,\IQ_p,\Rp:p\in\Pi\}$, where $\Pi$ is the set of prime
numbers, such that
\begin{itemize}
\item $\IQ\in\sigma(G)$ iff $G/\Tor(G)$ is divisible;
\item $\IZ_p\in\sigma(G)$ iff the $p$-torsion part $\pTor(G)$ is not divisible by $p$;
\item $\IQ_p\in\sigma(G)$ iff $\pTor(G)\ne0$ is divisible by $p$;
\item $\Rp\in\sigma(G)$ iff the group $G/\pTor(G)$ is not divisible by $p$.
\end{itemize}
Here $$\Tor(G)=\{x\in G:\exists n\in\IN\; n\cdot x=0\}\mbox{ and }\pTor(G)=\{x\in G:\exists k\in\IZ\; p^k\cdot x=0\}$$
is the torsion and $p$-torsion parts of $G$.
In particular, $\sigma(\IZ)=\{\Rp:p\in\Pi\}$.

\begin{theorem}\label{bockstein} Let $A$ be a closed subset of a space $X$, $G$ be a coefficient group,
and $p$ be a prime number.
\begin{enumerate}
\item $A$ is a $G$-homological $Z_n$-set in $X$ if and only if $A$ is an $H$-homological $Z_n$-set
in $X$ for all groups $H\in\sigma(G)$.
\item If $A$ is a $\Rp$-homological $Z_n$-set in $X$, then $A$ is a $\IQ$-homological
and $\IZ_p$-homological $Z_n$-set in $X$.
\item If $A$ is a $\IZ_p$-homological $Z_n$-set in $X$, then $A$ is a
$\IQ_p$-homological $Z_n$-set in $X$.
\item If $A$ is a $\IQ_p$-homological $Z_{n+1}$-set in $X$, then $A$ is a
$\IZ_p$-homological $Z_n$-set in $X$.
\item $A$ is a $\Rp$-homological $Z_n$-set in $X$ provided $A$ is a $\IQ$-homological
$Z_n$-set in $X$ and a $\IQ_p$-homological $Z_{n+1}$-set in $X$.
\end{enumerate}
\end{theorem}

By analogy with multiplication formulas for the
$\HDD[m]{n,k}$-properties, there are multiplication formulas for homotopical
and homological $Z_n$-sets, see \cite[Theorem
6.1]{BKV}).

\begin{theorem}\label{multiZZ} Let $A\subset X$, $B\subset Y$ be
closed subsets in Tychonoff spaces $X,Y$.
\begin{enumerate}
\item If $A$ is a homotopical $Z_n$-set in $X$ and $B$ is a
homotopical $Z_m$-set in $X$ then $A\times B$ is a homotopical
$Z_{n+m+1}$-set in $X\times Y$;
\item If $A$ is a homological
$Z_n$-set in $X$ and $B$ is a homological $Z_m$-set in $X$ then
$A\times B$ is a homological $Z_{n+m+1}$-set in $X\times Y$.
\end{enumerate}
\end{theorem}

Surprisingly, the multiplication formulas for homological $Z_n$-sets
can be reversed:

\begin{theorem}\label{multiZ} Let
$n,m\in\w\cup\{\infty\}$, $k\in\w$, and $A\subset X$, $B\subset Y$
be closed subsets of the spaces $X$ and $Y$. Let $\mathfrak
D=\{\IQ,\IQ_p:p\in\Pi\}$ and for every group
$G\in\mathfrak D$ let $B_G\subset Y$ be a closed subset which fails
to be a $G$-homological $Z_m$-set in $Y$. Then we have:
\begin{enumerate}
\item $A$ is a homological $Z_n$-set in $X$ if and only if $A^k$ is a homological $Z_{kn+k-1}$-set in
$X^k$;
\item If $A\times B$ is an $\IF$-homological $Z_{n+m}$-set in $X\times Y$for
some field $\IF$, then either $A$ is
an $\IF$-homological $Z_n$-set in $X$ or $B$ is an $\IF$-homological
$Z_m$-set in $Y$;
\item If $A\times B$ is a homological
$Z_{n+m}$-set in $X\times Y$, then either $A$ is a homological
$Z_n$-set in $X$ or $B$ is an $\exists G$-homological $Z_m$-set in
$Y$;
\item $A$ is a homological $Z_n$-set in $X$ provided $A\times B_G$ is a
homological $Z_{n+m}$-set in $X\times Y$ for every group
$G\in\mathfrak D$.
\end{enumerate}
\end{theorem}

Theorem~\ref{multiZ} is the principal tool for the proof of $k$-root
and multiplication formulas about the classes $\HDD[m]{n,k}$. We
first discuss $k$-root and division formulas for the classes $\Z_n$,
$\DZ_n$, and $\Z_n^{\IZ}$ because they are tightly connected with
the classes $\HDD[m]{n,k}$.

Let us start with some definitions. A point $x$ of a space $X$ is
defined to be a {\em homological $Z_n$-point\/} if its singleton
$\{x\}$ is a homological $Z_n$-set in $X$. By analogy, we define
$G$-homological and $\exists G$-homological $Z_n$-points.

Let $\Z_n^G(X)$ denote the set of all $G$-homological $Z_n$-points
in a space $X$ and $\DZ_n^G$ (resp., $\Z_n^G$) be the class of
Tychonoff spaces $X$ such that the set $\Z_n^G(X)$ is dense in
(resp., coincides with) $X$. We also recall that $\DZ_n$ (resp.,
$\mathcal Z_n$) stands for the class of Tychonoff spaces $X$ such
that the set $\Z_n(X)$ of homotopical $Z_n$-points of $X$ is dense
in (resp., coincides with) $X$. Using these notations,
Theorem~\ref{Z2+hZn=Zn} can be written in the following form.

\begin{theorem}\label{Zn-classes} Let $n\in\w\cup\{\infty\}$ and $G$ be a non-trivial Abelian group.
\begin{enumerate}
\item $\Z_n\subset\Z_n^{\IZ}\subset\Z_n^G$;
\item $\Z_0=\Z_0^{\IZ}=\Z_0^G$;
\item $\LC[1]\cap \Z_1^{G}\subset\Z_1$;
\item $\LC[1]\cap\Z_2\cap\Z_n^{\IZ}\subset\Z_n$;
\item $\LC[1]\cap\Br\cap \DZ_2\cap\DZ_n^{\IZ}\subset\DZ_n$.
\end{enumerate}
\end{theorem}

The last item of Theorem~\ref{Zn-classes} follows from the fact that
each of the sets $\Z_n(X)$ and $\Z_n^{\IZ}(X)$ is $G_\delta$ in $X$
provided $X$ is a separable metrizable $\LC[n]$-space \cite[Theorem
9.2]{BKV}.

In its turn, Theorem~\ref{multiZZ} implies multiplication formulas
for the classes $\Z_n$, $\DZ_n$, and $\Z_n^{\IZ}$:

\begin{theorem}[{\bf Multiplication Formulas}]\label{Zn-mult} Let $n,m\in\w\cup\{\infty\}$. Then
\begin{enumerate}
\item $\Z_n\times\Z_m\subset\Z_{m+n+1}$;
\item $\Z_n^{\IZ}\times\Z_m^{\IZ}\subset\Z^{\IZ}_{n+m+1}$;
\item $\DZ_n\times\DZ_m\subset\DZ_{m+n+1}$;
\item $\DZ_n^{\IZ}\times\DZ_m^{\IZ}\subset\DZ^{\IZ}_{n+m+1}$.

\end{enumerate}
\end{theorem}

The multiplication formulas can be reversed, which yields division
and $k$-root formulas for the classes $\Z_n^{\IZ}$ (we recall that,
for a class $\mathcal A$, we put $\sqrt[k]{\mathcal
A}=\{X:X^k\in\mathcal A\}$).

\begin{theorem}[{\bf $k$-Root Formulas}]\label{rootZ} Let $n\in\w\cup\{\infty\}$ and $k\in\IN$.
\begin{enumerate}
\item A space $X$ belongs to the class $\Z^{\IZ}_n$ if and only if $X^k$ belongs to $\Z_{kn+k-1}^{\IZ}$:
\smallskip

\hskip-10pt \frame{\phantom{$\Big|^{\big|}$}
$\Z_n^{\IZ}=\sqrt[k]{\Z^{\IZ}_{kn+k-1}}$
\phantom{$\Big|_{\big|}$}}
\smallskip

\item A metrizable separable Baire $\LC[kn+k-1]$-space $X$ belongs to the class
$\DZ_n^{\IZ}$ if and only if $X^k$ belongs to
$\DZ_{kn+k-1}^{\IZ}$\textup{:}
\smallskip

\hskip-10pt\frame{\phantom{$\Big|^{\big|}$}
$\DZ_n^{\IZ}\supset\sqrt[k]{\; \DZ^{\IZ}_{kn+k-1}}\cap \LC[nk+k-1]\cap\Br$
\phantom{$\Big|_{\big|}$}}
\smallskip
\end{enumerate}
\end{theorem}

To state the division formula for the classes $\Z^{\IZ}_n$ and
$\DZ^{\IZ}_n$ we need some more notations (which will be used  for
the classes $\HDD[m]{n,k}$ as well). Consider the following classes
of topological spaces:
\begin{itemize}
\item $\cup_G\Z^G_n=\cup\{\Z_n^G:G$ is a non-trivial Abelian group$\}$;
\item $\cup_G\DZ^G_n=\cup\{\DZ_n^G:G$ is a non-trivial Abelian group$\}$;
\item $\exists_\circ\kern-2pt\cup_G\kern-2pt\DZ^G_n$, the class of
spaces containing a non-empty open subspace $U\in\cup_G\DZ^G_n$.
\end{itemize}
For example, the space $\IR^n$ belongs to none of these classes.


Now, we can state the division
formulas for the classes $\Z^{\IZ}_n$ and $\DZ^{\IZ}_n$ (recall that
if $\mathcal A$ and $\mathcal B$ are two classes, then
$\displaystyle\frac{\mathcal A}{\mathcal B}$ stands for the class
$\{X\in\Top:\exists B\in\mathcal B\mbox{ with } X\times B\in\mathcal
A\}$).

\begin{theorem}[{\bf Division Formulas}]\label{Zn-div} Let $n\in\w\cup\{\infty\}$ and $k\in\w$.
\begin{enumerate}
\item A space $X$ belongs to the class $\Z^{\IZ}_n$ if and only if $X\times Y\in\Z^{\IZ}_{n+k}$
for some space $Y\notin\cup_G\Z^G_k$. This can be written as

\smallskip

\hskip-10pt\frame{\phantom{$\Big|^{\Big|}$}
 $\dfrac{\Z_{n+k}^{\IZ}}{\Top\setminus\cup_G\Z^G_k}= \Z^{\IZ}_n$
\phantom{$\Big|_{\Big|}$}}
\smallskip

\item A metrizable separable Baire $\LC[n]$-space $X$ belongs to the
class $\DZ^{\IZ}_n$ if and only if $X\times Y\in\DZ^{\IZ}_{n+k}$ for
some space $Y\notin\exists_\circ\kern-2pt\cup_G\DZ^G_n$.
This can be written as
\smallskip

\hskip-10pt\frame{\phantom{$\Big|^{\Big|}$} $\Br\cap\LC[n]\cap\dfrac{\DZ_{n+k}^{\IZ}}
{\Top\setminus\exists_\circ\kern-2pt\cup_G\kern-2pt\DZ^G_k}\subset \DZ^{\IZ}_n$ \phantom{$\Bigg|_{\big|}$}}
\smallskip

\end{enumerate}
\end{theorem}

Because of the division formulas, it is important to detect the
spaces $X\notin\cup_G \Z^G_n$. It turns out that this happens for
every metrizable space $X$ of dimension $\dim X\le n$, or more
generally of transfinite separation dimension $\trt(X)<n+1$. The
latter dimension can be introduced  inductively (see \cite{ACP}):
\begin{itemize}
\item $\trt(X)=-1$ iff $X=\varnothing$;
\item $\trt(X)\le\alpha$ for an ordinal $\alpha$ if each closed subset $A\subset X$ with
$|A|>1$ contains a closed subset $B\subset A$ such that $\trt(B)<\alpha$ and $A\setminus B$ is disconnected.
\end{itemize}
A space $X$ is called {\em $\trt$-dimensional} if $\trt(X)\le\alpha$
for some ordinal $\alpha$. For a $\trt$-dimensional space $X$ we put
$\trt(X)$ be the smallest ordinal $\alpha$ with $\trt(X)\le\alpha$.

By \cite{ACP}, each compact metrizable $\trt$-dimensional space is a
$C$-space. On the other hand, a \v Cech-complete space is
$\trt$-dimensional if it can be written as the countable union of
hereditarily disconnected subspaces, see \cite{Ra}. It
is easy to see that for a finite-dimensional metrizable separable
space $X$ we get $\trt(X)\le\dim(X)$. Moreover, if $X$ is
finite-dimensional and compact, then $\trt(X)=\dim(X)$, see
\cite{St}.

The following theorem was proved in \cite{BKV} and \cite{BC}.

\begin{theorem}\label{dimZn} Let $X\in\cup_G\Z^G_n$ for some $n\in\w\cup\{\infty\}$.
\begin{enumerate}
\item If $n<\infty$, then $\trt(X)>n$;
\item If $n=\infty$, then $X$ is not $\trt$-dimensional;
\item If $n=\infty$ and $X$ is locally compact and locally contractible, then $X$ is not a $C$-space.
\end{enumerate}
\end{theorem}

Consequently, for any metrizable separable space $X\in\cup_G\Z_n^G$ we have
$\dim X\ge\trt(X)>n$. A similar inequality holds for cohomological and
extension dimensions of $X$. We recall their definitions.

 For a  space $X$ and a $CW$-complex $L$ we
write $\edim X\le L$ if each map $f:A\to L$ defined on a closed
subset $A\subset X$ admits a continuous extension $\bar f:X\to L$,
see \cite{DD} for more information on Extension Dimension Theory. It
follows from the classical Hurewicz-Wallman Theorem
\cite[1.9.3]{End} that $\edim X \le S^n$ iff $\dim X\le n$. The
cohomological dimension with respect to a given Abelian group $G$
can be expressed via extension dimension as follows: define
$\dim_GX\le n$ if $\edim X\le K(G,n)$, where $K(G,n)$ is the
Eilenberg-MacLane complex of type $(G,n)$, and let $\dim_GX$ be the
smallest non-negative integer with $\dim_GX\le n$. If there is no
such an integer $n$, we put $\dim_GX=\infty$.

\begin{theorem}\label{edimZn} Let $n\in\w$ and $X\in \Z^{\IZ}_n$ be
a locally compact $\LC[n]$-space. Then
\begin{enumerate}
\item $\dim_GX> n$ for any Abelian group $G$;
\item $\edim X\not\le L$ for any $CW$-complex $L$ with a
non-trivial homotopy group $\pi_k(L)$ for some $k\le n$.
\end{enumerate}
\end{theorem}

\subsection{Homological $Z_n$-sets and $\HDD[m]{n,k}$-properties}
In this subsection we discuss the interplay between the classes
$\Z_n^G$ and $\HDD[m]{n,k}$.  The following two theorems present
homological counterparts of the formulas
$$\begin{aligned}
&\HDD[n]{n,0}\;\subset\; \Z_n\;\subset\;\bigcap_{m+k\le
n}\HDD[m]{0,k}\mbox{ and}\\
&\Br\cap \LC[n]\cap \HDD[0]{0,n}\;\subset\;
\LC[0]\cap\DZ_n\;\subset\;
\HDD[0]{0,n}
\end{aligned}$$
from Theorem~\ref{Zset}.

\begin{theorem}\label{HDD->Z} Let $X$ be a Tychonoff space and
$n\in\w\cup\{\infty\}$.
\begin{enumerate}
\item If an $\LC[1]$-space $X$ has the $\HDD[2]{0,2}$-property,
then each homological $Z_n$-point in $X$ is a homotopical
$Z_n$-point:
\smallskip

\hskip-15pt\frame{\phantom{$I^{I^{I^{I^I}}}$}\hskip-15pt
$\LC[1]\cap\HDD[2]{0,2}\cap\Z_n^{\IZ}\subset\Z_n$
\hskip-15pt\phantom{$I_{I_{I_{I_I}}}$}}
\smallskip

\item If a metrizable separable Baire $\LC[n]$-space $X$ has the $\HDD[0]{0,2}$-property
and contains a dense set of homological $Z_n$-points, then $X$
contains a dense set of homotopical $Z_n$-points and
$X\in\HDD[0]{0,n}$:
\smallskip

\hskip-15pt\frame{\phantom{$I^{I^{I^{I^I}}}$}\hskip-15pt
$\LC[n]\cap\Br\cap\HDD[0]{0,2}\cap\DZ_n^{\IZ}\subset\DZ_n$
\hskip-15pt\phantom{$I_{I_{I_{I_I}}}$}}
\smallskip

\item If $X$ has the $\HDD[(2n+1)]{0,0}$-property, then each point
of $X$ is a homological $Z_n$-point:
\smallskip

\hskip-15pt\frame{\phantom{$I^{I^{I^{I^I}}}$}\hskip-15pt
$\HDD[(2n+1)]{0,0}\subset \Z_n^{\IZ}$
\hskip-15pt\phantom{$I_{I_{I_{I_I}}}$}}
\smallskip

\item If $X$ has the $\HDD[2n]{0,0}$-property, then each point of
$X$ is a $G$-homological $Z_n$-point for any group $G$ with
divisible quotient $G/\Tor(G)$. Consequently,
\smallskip

\hskip-15pt\frame{\phantom{$I^{I^{I^{I^I}}}$}\hskip-15pt
$\HDD[2n]{0,0}\subset \Z_n^{\IQ}\cap
\bigcap_p(\Z^{\IZ_{p}}_n\cap\Z^{\IQ_p}_n)$
\hskip-15pt\phantom{$I_{I_{I_{I_I}}}$}}
\smallskip
\end{enumerate}
\end{theorem}

\begin{theorem}\label{Z->HDD} Let $m,m,k$ be non-negative integers or infinity.
\begin{enumerate}
\item If each point of an $\LC[1]$-space $X$ with the
$\HDD[2]{0,2}$-property is a homological $Z_{m+k}$-point, then $X$
has the $\HDD[m]{0,k}$-property.  This can be written as
\smallskip

\hskip-15pt\frame{\phantom{$I^{I^{I^{I^I}}}$}\hskip-15pt
$\LC[1]\cap\Z_{m+k}^{\IZ}\cap\HDD[2]{0,2}
\subset\HDD[m]{0,k}$\hskip-15pt\phantom{$I_{I_{I_{I_I}}}$}}
\smallskip

\item If a metrizable separable $\LC[k]$-space $X$ has the
$\HDD[0]{0,2}$-property and contains a dense set of homological
$Z_{k}$-points, then $X$ has the $\HDD[0]{0,k}$-property. This can
be written as
\smallskip

\hskip-15pt\frame{\phantom{$I^{I^{I^{I^I}}}$}\hskip-15pt
$\LC[k]\cap\DZ_{k}^{\IZ}\cap \HDD[0]{0,2}\subset
\HDD[0]{0,k}$\hskip-15pt\phantom{$I_{I_{I_{I_I}}}$}}
\smallskip

\item If each point of a metrizable separable $\LC[2]$-space $X$ is
a homological $Z_{m+n+k}$-point and $X$ has the
$\HDD[m]{n,\max\{n,2\}}$-property, then $X$ has the
$\HDD[m]{n,k}$-property. This can be written as
\smallskip

\hskip-15pt\frame{\phantom{$I^{I^{I^{I^I}}}$}\hskip-15pt
$\LC[2]\cap\Z_{m+n+k}^{\IZ}\cap \HDD[m]{n,\max\{2,n\}}\subset
\HDD[m]{n,k}$\hskip-15pt\phantom{$I_{I_{I_{I_I}}}$}}
\smallskip

\item If each point of an  $\LC[1]$-space $X$ is a homological
$Z_m$-point and $X\in \HDD[2]{0,0}$, then $X$ has the
$\HDD[m]{0,0}$-property. This can be written as
\smallskip

\hskip-15pt\frame{\phantom{$I^{I^{I^{I^I}}}$}\hskip-15pt
$\LC[1]\cap\Z_{m}^{\IZ}\cap \HDD[2]{0,0}\subset
\HDD[m]{0,0}$\hskip-15pt\phantom{$I_{I_{I_{I_I}}}$}}
\smallskip

\item If each point of an $\LC[0]$-space $X$ is a $G$-homological
$Z_2$-point for some group $G$, then $X$ has the $\mathrm{DD^1P}$.
If in addition $X$ is a metrizable $\LC[1]$-space containing a dense
set of homotopical $Z_2$-points, then $X\in\HDD[2]{0,0}$.
\smallskip

\hskip-15pt\frame{\phantom{$I^{I^{I^{I^I}}}$}\hskip-16pt $\LC[0]\cap
(\cup_G\Z_{2}^{G})\subset \HDD[0]{1,1}$ \;\; and \;\;
$\LC[1]\cap\DZ_2\cap (\cup_G\Z_2^G)\subset \HDD[2]{0,0}$
\hskip-20pt\phantom{$I_{I_{I_{I_I}}}$}}
\end{enumerate}
\end{theorem}

\subsection{Homological characterization of the $\HDD[0]{n,k}$-property}
In this subsection we prove a quantified version of the homological
characterization of the $\HDD[0]{\infty,\infty}$-property due to
Daverman and Walsh \cite{DW}.

First, we provide a homotopical version of the Daverman-Walsh
result.

\begin{theorem}\label{h-0DDnk}
Let $n,k$ be finite or infinite integers \textup{(}with $n\le
k$\textup{)}. A Polish \textup{(}$\LC[k]$\textup{)}-space $X$ has
the $\HDD[0]{n,k}$-property if \textup{(}and only if\textup{)} there
is a countable family $\F$ of \textup{(}$n$-dimensional
compact\textup{)} homotopical $Z_k$-sets in $X$ such that each
compact subset $K\subset X\setminus\cup\F$ is a homotopical
$Z_{n}$-set in $X$.
\end{theorem}

Under some mild assumptions on $X$ it is possible to replace the
homotopical conditions in Theorem~\ref{h-0DDnk} by homological ones.

\begin{theorem}\label{DW-0DDnk}  A Polish $\LC[\max\{n,k\}]$-space
$X\in\HDD[0]{2,2}$ has the $\HDD[0]{n,k}$-property provided each
point of $X$ is a homological $Z_{2+\max\{n,k\}}$-point and there is
a countable family $\F$ of homological $Z_k$-sets in $X$ such that
each compact subset $K\subset X\setminus\cup\F$ is a homological
$Z_n$-set in $X$.
\end{theorem}

This theorem implies another characterization of $\HDD[0]{n,k}$ in
terms of approximation properties defined as follows. We shall say
that a topological space has the {\em $n$-dimensional approximation
property} (briefly, $\AP[n]$) if for any open cover $\U$ of $X$ and
a map $f:\mathbb I^n\to X$ there is a map $g:\mathbb I^n\to X$ such
that $g$ is $\U$-homotopic to $f$ and $\trt(g(\mathbb I^n))<n+1$.
Here we assume that $\alpha<\infty+1$ for each ordinal $\alpha$
(which is essential if $n=\infty$).

Observe that each $\LC[0]$-space has $\AP[0]$ and each $\LC[1]$-space has $\AP[1]$.

\begin{theorem}\label{AP->HDD}
If each point of a Polish $\LC[\max\{n,k\}]$-space $X$ is a
homological $Z_{n+k}$-point and $X$ has the properties $\AP[n]$ and
$\HDD[0]{2,\min\{2,n\}}$, then $X$ has the $\HDD[0]{n,k}$-property.
This can be written as
\smallskip

\hskip10pt\frame{\phantom{$I^{I^{I^{I^I}}}$}\hskip-15pt
$\Pi^0_2\cap \LC[\max\{n,k\}]\cap\Z_{n+k}^{\IZ}\cap
\AP[n]\cap\HDD[0]{2,\min\{2,n\}}\subset
\HDD[0]{n,k}$\hskip-15pt\phantom{$I_{I_{I_{I_I}}}$}}
\smallskip
\end{theorem}

\subsection{$\HDD[m]{n,k}$-properties of locally rectifiable spaces}

There is a non-trivial interplay between $\HDD[m]{n,k}$-properties
for spaces having a kind of a homogeneity property. We recall that a
space $X$ is {\em topologically homogeneous} if for any two points
$x_0,x\in X$ there is a homeomorphism $h_{x}:X\to X$ such that
$h_{x}(x_0)=x$. If the homeomorphism $h_x$ can be chosen to depend
continuously on $x$ then $X$ is called rectifiable at $x_0$.

More precisely, we define a topological space $X$ to be {\em locally
rectifiable at a point $x_0\in X$} if there exists a neighborhood
$U$ of $x_0$ such that for every $x\in U$ there is a homeomorphism
$h_x:X\to X$ such that $h_x(x_0)=x$ and $h_x$ continuously depends
on $x$ in the sense that the map $H:U\times X\to U\times X$,
$H:(x,z)\mapsto (x,h_x(z))$ is a homeomorphism. If $U=X$, then the
space $X$ is called {\em rectifiable} at $x_0$.

A space $X$ is called ({\em locally}) {\em rectifiable} if it is
(locally) rectifiable at each point $x\in X$. Rectifiable spaces
were studied in details by A.S.~Gulko \cite{Gul}. The class of
rectifiable spaces contains the underlying spaces of topological
groups but also contains spaces not homeomorphic to topological
groups. A simplest such an example is the 7-dimensional sphere
$S^7$, see \cite{Us89}. It should be mentioned that all
finite-dimensional spheres $S^n$ are locally rectifiable but only
$S^1,S^3$ and $S^7$ are rectifiable (this follows from the famous
Adams' result \cite{Adams} detecting $H$-spaces among the spheres).
It can be shown that each connected locally rectifiable space is
topologically homogeneous. On the other hand, the Hilbert cube is
topologically homogeneous but fails to be (locally) rectifiable, see
\cite{Gul}.

By $\lCH$ we denote the class of Tychonoff locally rectifiable spaces.

\begin{theorem}\label{homogen} Let $X$  be a locally rectifiable Tychonoff space.
\begin{enumerate}
\item If $X$ has the $\HDD[m]{0,k}$-property, then each point of
$X$ is a homotopical $Z_{m+k}$-point:
\smallskip

\hskip-10pt\frame{\phantom{$I^{I^{I^{I^I}}}$}\hskip-15pt $\lCH\cap
\HDD[m]{0,k}\subset \Z_{m+k}$
\hskip-15pt\phantom{$I_{I_{I_{I_I}}}$}}
\smallskip

\item If $X$ has the $\HDD[m]{0,k}$-property, then $X$ has
$\HDD[i]{0,j}$-properties for all $i,j$ with $i+j\le m+k$:
\smallskip

\hskip-10pt\frame{\phantom{$I^{I^{I^{I^I}}}$}\hskip-15pt $\lCH\cap
\HDD[m]{0,k}\subset \bigcap_{i+j\le m+k}\HDD[i]{0,j}$
\hskip-15pt\phantom{$I_{I_{I_{I_I}}}$}}
\smallskip

\item If either $X\in\DZ_{m+p}$ or $X\in
\DZ_{m+p}^{\IZ}\cap\LC[1]$, then the product $X\times Y$ has the
$\HDD[m]{n,k+p+1}$-property for each separable metrizable
$\LC[k]$-space $Y$ possessing the $\HDD[m]{n,k}$-property with $n\le
k$. This can be written as
\smallskip

\hskip-10pt\frame{\phantom{$I^{I^{I^{I^I}}}$} \hskip-15pt
$\big(\lCH\cap \LC[1]\cap \DZ^{\IZ}_{m+p}\big) \times
\big(\LC[k]\cap\HDD[m]{n,k}\big)\subset\HDD[m]{n,k+p+1}$
\hskip-15pt\phantom{$I_{I_{I_{I_I}}}$}}
\end{enumerate}
\end{theorem}

\begin{remark}
Since $\mathbb R^q\in \Z_{q-1}$ is rectifiable,
Theorem~\ref{homogen}(3) implies that the product $X\times\mathbb
R^{m+p}$ has the $\HDD[m]{n,k+p}$-property for any separable
metrizable $\LC[k]$-space having the $\HDD[m]{n,k}$-property with
$n\leq k$. This result was established by W.~Mitchell \cite[Theorem
4.3(1)]{mitchell} in the case $X$ is a compact $\ANR$ and $m=0$.
Moreover, a particular case of Theorem~\ref{homogen}(1) when $X$ is
an $\ANR$ and $m=0$ was also established in \cite{mitchell}.
\end{remark}

\subsection{$k$-Root and Division Formulas for the $\HDD[m]{n,k}$-properties}
In this section we discuss $k$-root and division formulas for the
$\HDD[m]{n,k}$-properties, one of the most surprising features of
these properties.

\begin{theorem}\label{root} {\bf ($k$-Root Formulas)}
Let $n$ be a non-negative integer or infinity and $k$ be a positive integer.
\begin{enumerate}
\item
If $X$ is an $\LC[1]$-space with $\HDD[2]{0,0}$-property and
$X^k\in\HDD[(kn{+}k{-}1)]{0,0}$, then $X$ has the $\HDD[n]{0,0}$-property.
 This can be written as
\medskip

\hskip-10pt\frame{\phantom{\huge $A^{A^a}$}\hskip-20pt
$\LC[1]\cap
\HDD[2]{0,0}\cap\sqrt[k]{\HDD[{\scriptsize{(kn+k-1)}}]{0,0}}\subset\HDD[n]{0,0}$
\hskip-20pt\phantom{\huge $A_{|_|}$}}

\item
If $X$ is a separable metrizable $\LC[kn+k-1]$-space
with the $\HDD[0]{0,2}$-property and $X^k$ has the
$\HDD[0]{0,kn+k-1}$-property, then $X$ has the $\HDD[0]{0,n}$-property.
 This can be written as
\medskip

\hskip-10pt\frame{\phantom{\huge $A^{A^a}$}\hskip-15pt
$\LC[kn+k-1]\cap\HDD[0]{0,2}\cap\sqrt[k]{\phantom{I^{I^I}}\kern-15pt\HDD[0]{0,kn+k-1}}
\subset\HDD[0]{0,n}$
\hskip-15pt\phantom{\huge $A_{|_l}$}}
\end{enumerate}
\end{theorem}

To write down division formulas for the $\HDD[m]{n,k}$-property,
let us introduce two new classes in addition to the classes
$\cup_G\Z^G_n$ and $\exists_\circ\kern-2pt \cup_G\DZ^G_n$:
\begin{itemize}
\item $\Z^{\exists G}_n$ - the class of spaces $X$ with all $x\in X$
being $\exists G$-homological $Z_n$-points in $X$;
\item $\Delta\Z_n^{\exists G}$ - the class of spaces $X$ whose diagonal
$\Delta_X$ is an $\exists G$-homological $Z_n$-set in $X^2$.
\end{itemize}
Note that any at most $n$-dimensional polyhedron belongs to none of
the last two classes.

\begin{theorem}\label{division} {\bf (Division Formulas)}
Let $n\le k$ be non-negative integers or infinity and $m$ a non-negative integer.
\begin{enumerate}
\item
An $\LC[1]$-space with the $\HDD[2]{0,0}$-property has the
$\HDD[n]{0,0}$-property provided $X\times Y$ has the
$\HDD[(n+m)]{0,0}$-property for some space $Y$ whose diagonal
$\Delta_Y$ fails to be a $\exists G$-homological $Z_m$-set in $Y^2$.
This can be written as
\medskip

\hskip-15pt\frame{\phantom{$I^{I^{I^{I^{I^{I^{I^i}}}}}}$}\hskip-20pt
$\LC[1]\cap \HDD[2]{0,0}\cap\dfrac{\HDD[(n+m)]{0,0}}{\Top\setminus
\Delta\Z^{\exists G}_m}\subset\HDD[n]{0,0}$
\hskip-20pt\phantom{$I_{I_{I_{I_{I_{I_{I_{I_I}}}}}}}$}}

\item
A separable metrizable $\LC[n+m]$-space $X\in\HDD[0]{0,2}$ has the
$\HDD[0]{0,n}$-property provided $X\times Y$ has the
$\HDD[0]{0,n+m}$-property for some metrizable separable Baire
$\LC[n+m]$-space $Y$ that contains no non-empty open set $U\in
\cup_G\DZ^G_m$.
\medskip

\hskip-15pt\frame{\phantom{$I^{I^{I^{I^{I^{I^{I^i}}}}}}$}\hskip-20pt
$\LC[n+m]\cap\HDD[0]{0,2}\cap\dfrac{\HDD[0]{0,n+m}}
{\Br\cap\LC[n+m]\setminus\exists_\circ\kern-2pt\cup_G\kern-2pt\DZ^G_m}\subset\HDD[0]{0,n}$
\hskip-10pt\phantom{\huge $|_{|_{|_{|_|}}}$}}

\item A separable metrizable $\LC[k+m]$-space $X\in \Z^{\IZ}_{k+2}$
with the $\HDD[0]{2,2}$-property has the $\HDD[0]{n,k}$-property
provided $X\times Y$ has the $\HDD[0]{n+m,k+m}$-property for some
metrizable separable $\LC[k+m]$-space $Y\notin\Z^{\exists
G}_m$. This can be written as
\medskip

\hskip-15pt\frame{\phantom{$I^{I^{I^{I^{I^{I^{I^i}}}}}}$}\hskip-20pt
$\LC[k+m]\cap\HDD[0]{2,2}\cap\Z^{\IZ}_{k+2}\cap\dfrac{\HDD[0]{n+m,k+m}}
{\LC[k+m]\setminus\Z_m^{\exists
G}}\subset\HDD[0]{n,k}$
\hskip-20pt\phantom{$I_{I_{I_{I_{I_{I_{I_{I_I}}}}}}}$}}

\item A separable metrizable $\LC[k+m]$-space
$X\in\Z^{\IZ}_{n+k+m}$ possessing the $\HDD[0]{2,2}$-property has
the $\HDD[0]{n,k}$-property provided $X\times Y$ has the
$\HDD[0]{n+m,n+m}$-property for some metrizable separable
$\LC[n+m]$-space $Y\notin\cup_G\Z^G_m$. This can be written as
\medskip

\hskip-15pt\frame{\phantom{$I^{I^{I^{I^{I^{I^{I^i}}}}}}$}\hskip-20pt
$\LC[k+m]\cap\HDD[0]{2,2}\cap\Z^{\IZ}_{n+k+m}\cap\dfrac{\HDD[0]{n+m,n+m}}
{\LC[n+m]\setminus\cup_G\Z^G_m}\subset\HDD[0]{n,k}$
\hskip-20pt\phantom{$I_{I_{I_{I_{I_{I_{I_{I_I}}}}}}}$}}
\end{enumerate}
\end{theorem}

\subsection{Characterizing $\HDD[m]{n,k}$-properties with $m,n,k\in\{0,\infty\}$}
In this subsection we apply the results obtained in preceding
subsections to the case of $\HDD[m]{n,k}$-properties with
$m,n,k\in\{0,\infty\}$. Let us note that the $\HDD[0]{0,0}$ has been
characterized in Proposition~\ref{lc-00} while
$\HDD[\infty]{\infty,\infty}$ is equivalent to
$\HDD[0]{\infty,\infty}$. So, it suffices to consider only the
properties: $\HDD[0]{0,\infty}$, $\HDD[\infty]{0,0}$,
$\HDD[\infty]{0,\infty}$, and $\HDD[0]{\infty,\infty}$. These
properties can be characterized in terms of homotopical or
homological $Z_\infty$-points as follows:

\begin{corollary}\label{infty-HDD} \phantom{mm}

\begin{enumerate}
\item A topological \textup{(}$\LC[1]$-\textup{)}space $X$ has the
$\HDD[\infty]{0,\infty}$-property \newline iff all points of $X$ are
homotopical $Z_\infty$-points \newline\textup{(}iff all points of
$X$ are homotopical $Z_2$-points and $X$ has
$\HDD[\infty]{0,0}$-property\textup{)}:

\smallskip

\hskip-10pt\frame{\phantom{$I^{I^{I^{I^I}}}$}\hskip-15pt
$\Z_2\cap\HDD[\infty]{0,0}\cap\LC[1]\subset
\HDD[\infty]{0,\infty}=\Z_\infty$
\hskip-15pt\phantom{$I_{I_{I_{I_I}}}$}}
\smallskip

\item An $\LC[1]$-space $X$ has the $\HDD[\infty]{0,0}$-property
\newline  iff $X\in\HDD[2]{0,0}$ and all points of $X$ are
homological $Z_\infty$-points:
\smallskip

\hskip-10pt\frame{\phantom{$I^{I^{I^{I^I}}}$}\hskip-15pt
$\Z_\infty^{\IZ}\cap\HDD[2]{0,0}\cap\LC[1]\subset
\HDD[\infty]{0,0}\subset\Z^{\IZ}_\infty$
\hskip-15pt\phantom{$I_{I_{I_{I_I}}}$}}
\smallskip

\item A Polish $\LC[\infty]$-space $X$ has the
$\HDD[0]{0,\infty}$-property\newline  iff $X$ has a dense set of
homotopical $Z_\infty$-points\newline  iff $X\in\HDD[0]{0,2}$ and
$X$ has a dense set of homological $Z_\infty$-points:
\smallskip

\hskip-25pt\frame{\phantom{$I^{I^{I^{I^I}}}$}\hskip-20pt
$\DZ^{\IZ}_\infty\cap\HDD[0]{0,2}\cap\LC[\infty]\subset
\HDD[0]{0,\infty}$\ and \ $\HDD[0]{0,\infty}\cap
\LC[\infty]\cap\Pi^0_2\subset
\DZ_\infty$\hskip-18pt\phantom{$I_{I_{I_{I_I}}}$}}

\item If each point of a metrizable separable $\LC[\infty]$-space
$X$ is a homological $Z_{\infty}$-point and $X$ has the properties
$\AP[\infty]$ and $\HDD[0]{2,2}$, then $X$ has the
$\HDD[0]{\infty,\infty}$-property:
\smallskip

\hskip-10pt\frame{\phantom{$I^{I^{I^{I^I}}}$}\hskip-15pt
$\Z_{\infty}^{\IZ}\cap\HDD[0]{2,2}\cap \LC[\infty]\cap\AP[\infty]\subset
\HDD[\infty]{\infty,\infty}=\HDD[0]{\infty,\infty}$\hskip-15pt\phantom{$I_{I_{I_{I_I}}}$}}
\end{enumerate}
\end{corollary}

According to the famous characterization of Hilbert cube manifolds
due to Toru\'nczyk \cite{To80}, a locally compact $\ANR$-space $X$
is a $Q$-manifold if and only if $X$ has the
$\HDD[0]{\infty,\infty}$-property. Combining this characterization
with the last item of Corollary~\ref{infty-HDD} we obtain a new
characterization of Q-manifolds.

\begin{corollary} A locally compact $ANR$ is a $Q$-manifold if and only if
\begin{itemize}
\item $X$ has the disjoint disk property;
\item each point of $X$ is a homological $Z_\infty$-point;
\item each map $f:\mathbb I^\infty\to X$ can be uniformly
approximated by maps with $\trt$-dimensional image.
\end{itemize}
\end{corollary}

Next, we discuss $k$-root formulas for $\HDD[m]{n,k}$-properties
with $m,n,k\in\{0,\infty\}$.

\begin{corollary}\label{infty-root} {\bf ($k$-Root Formulas)}
\begin{enumerate}
\item An $\LC[1]$-space $X$ has the $\HDD[\infty]{0,0}$-property if
and only if $X\in\HDD[2]{0,0}$ and $X^k\in\HDD[\infty]{0,0}$ for
some finite $k$. This can be written as
\smallskip

\hskip-10pt\frame{\phantom{$I^{I^{I^{I^I}}}$}\hskip-15pt
$\HDD[\infty]{0,0}\supset {\displaystyle
\sqrt[k]{\HDD[\infty]{0,0}}}\;\cap \HDD[2]{0,0}\cap \LC[1]$
\hskip-15pt\phantom{$I_{I_{I_{I_I}}}$}}
\smallskip

\item A metrizable separable $\LC[\infty]$-space $X$ has the
$\HDD[0]{0,\infty}$-property if and only if $X\in\HDD[0]{0,2}$ and
$X^k\in\HDD[0]{0,\infty}$ for some finite $k$. This can be written
as
\smallskip

\hskip-10pt\frame{\phantom{$I^{I^{I^{I^I}}}$}\hskip-15pt
$\HDD[0]{0,\infty}\supset{\displaystyle
\sqrt[k]{\HDD[0]{0,\infty}}}\;\cap\HDD[0]{0,2}\cap \LC[\infty]$
\hskip-15pt\phantom{$I_{I_{I_{I_I}}}$}}
\smallskip

\item An $\LC[1]$-space $X$ has the
$\HDD[\infty]{0,\infty}$-property iff $X$ has the
$\HDD[2]{0,2}$-property and $X^k\in\HDD[\infty]{0,\infty}$ for some
finite $k$. This can be written as
\smallskip

\hskip-10pt\frame{\phantom{$I^{I^{I^{I^I}}}$}\hskip-15pt
$\HDD[\infty]{0,\infty}\supset{\displaystyle
\sqrt[k]{\HDD[\infty]{0,\infty}}}\;\cap\HDD[2]{0,2}\cap \LC[1]$
\hskip-15pt\phantom{$I_{I_{I_{I_I}}}$}}
\smallskip

\item A metrizable separable $\LC[\infty]$-space $X$ has the
$\HDD[0]{\infty,\infty}$-property if $X$ has the
$\HDD[0]{2,2}$-property, $X\in\AP[\infty]$ and
$X^k\in\HDD[0]{\infty,\infty}$ for some finite $k$. This can be
written as
\smallskip

\hskip-10pt\frame{\phantom{$I^{I^{I^{I^I}}}$}\hskip-15pt
$\HDD[0]{\infty,\infty}\supset{\displaystyle \sqrt[k]{\HDD[0]{\infty,\infty}}}\;
\cap\HDD[0]{2,2}\cap \LC[\infty]\cap\AP[\infty]$ \hskip-15pt\phantom{$I_{I_{I_{I_I}}}$}}
\end{enumerate}
\end{corollary}

Finally, we turn to division formulas for the
$\HDD[m]{n,k}$-properties with $m,n,k\in\{0,\infty\}$.

\begin{corollary}\label{infty-div} {\bf (Division Formulas)}

\begin{enumerate}
\item An $\LC[1]$-space $X$ has the $\HDD[\infty]{0,0}$-property
provided $X$ has the $\HDD[2]{0,0}$-property and the product
$X\times Y$ has the $\HDD[\infty]{0,0}$-property for some space
$Y\notin \cup_G\Z^{G}_\infty$. This can be written as
\smallskip

\hskip-10pt\frame{\phantom{$I^{I^{I^{I^{I^{I^{I}}}}}}$}\hskip-15pt
$\HDD[\infty]{0,0}\supset \dfrac{\HDD[\infty]{0,0}}{\Top\setminus
\cup_G\Z^{G}_\infty}\cap \HDD[2]{0,0}\cap\LC[1]$
\hskip-15pt\phantom{$I_{I_{I_{I_{I_{I_{I_{I}}}}}}}$}}
\smallskip

\item An $\LC[1]$-space $X$ has the
$\HDD[\infty]{0,\infty}$-property provided $X$ has the
$\HDD[2]{0,2}$-property and the product $X\times Y$ has the
$\HDD[\infty]{0,\infty}$-property for some space $Y\notin
\cup_G\Z^{G}_\infty$. This can be written as
\smallskip

\hskip-10pt\frame{\phantom{$I^{I^{I^{I^{I^{I^{I}}}}}}$}\hskip-15pt
$\HDD[\infty]{0,\infty}\supset
\dfrac{\HDD[\infty]{0,\infty}}{\Top\setminus
\cup_G\Z^{G}_\infty}\cap \HDD[2]{0,2}\cap\LC[1]$
\hskip-15pt\phantom{$I_{I_{I_{I_{I_{I_{I_{I}}}}}}}$}}
\smallskip

\item A metrizable separable $\LC[\infty]$-space $X$ has the
$\HDD[0]{\infty,\infty}$-property provided $X\in\HDD[0]{2,2}$ and
$X\times Y$ has the $\HDD[0]{\infty,\infty}$-property for some
separable metrizable $\LC[\infty]$-space $Y\notin\cup_G\Z^G_\infty$.
This can be written as
\smallskip

\hskip-10pt\frame{\phantom{$I^{I^{I^{I^{I^{I^I}}}}}$}\hskip-15pt
$\HDD[0]{\infty,\infty}\supset\dfrac{\LC[\infty]\cap
\HDD[0]{\infty,\infty}}{\Top\setminus\cup_G\Z^G_\infty}\cap\HDD[0]{2,2}$
\hskip-15pt\phantom{$I_{I_{I_{I_{I_{I_{I_I}}}}}}$}}
\smallskip

\item A metrizable separable $\LC[\infty]$-space $X$ has the
$\HDD[0]{0,\infty}$-property provided $X\in\HDD[0]{0,2}$ and the
product $X\times Y$ has the $\HDD[0]{0,\infty}$-property for some
metrizable separable $\LC[\infty]$-space $Y\notin
\exists_\circ\kern-2pt\cup_G\kern-2pt\DZ^G_\infty$. This can be
written as
\smallskip

\hskip-10pt\frame{\phantom{$I^{I^{I^{I^I{^{I^{I^I}}}}}}$}\hskip-15pt
$\HDD[0]{0,\infty}\supset\dfrac{\LC[\infty]\cap\HDD[0]{0,\infty}}
{\Top\setminus \exists_\circ\kern-2pt\cup_G\kern-2pt\DZ^G_\infty}
\cap\HDD[0]{0,2}$
\hskip-15pt\phantom{$I_{I_{I_{I_{I_{I_{I_{I_{I_I}}}}}}}}$}}

\end{enumerate}
\end{corollary}

The third item of Corollary~\ref{infty-div} combined with the
Toru\'nczyk's characterization of $Q$-manifolds implies the
following division theorem for $Q$-manifolds proven in \cite{BC} and
implicitly in \cite{DW}.

\begin{corollary} A space $X$ is a $Q$-manifold if and only if
the product $X\times Y$ is a $Q$-manifold for some space
$Y\notin\cup_G\Z^G_\infty$.
\end{corollary}

\subsection{Dimension of spaces with the $\HDDP[m]{n}$-property}
In this section we study the dimensional properties of spaces
possessing the $\HDDP[m]{n}$-property.

\begin{theorem}\label{dim}
If a metrizable separable space $X$ has the $\HDDP[m]{n}$-property,
then $\displaystyle\dim X\ge n+\frac{m+1}{2}$.
\end{theorem}

This theorem combined with Theorem~\ref{power} allows us to
calculate the smallest possible dimension of a space $X$ with
$\HDDP[m]{n}$. For a real number $r$ let
$$\begin{aligned}
\lfloor r\rfloor=&\max\{n\in\IZ:n\le r\}\\
\lceil r\rceil=&\min\{n\in\IZ:n\ge r\}.
\end{aligned}
$$

\begin{corollary}\label{bound} Let $n,m$ be non-negative integers
and $D$ be a dendrite with a dense set of end-points.
\begin{enumerate}
\item If $m$ is odd and $\displaystyle d=n+\frac{m+1}{2}$, then the
power $D^d$ is a $d$-dimensional absolute retract with the
$\HDDP[m]{n}$-property.
\item If $m$ is even and $\displaystyle
d=n+\frac{m+2}2$, then the product $D^{d-1}\times\mathbb I$ is a
$d$-dimensional absolute retract with the $\HDDP[m]{n}$-property.
\end{enumerate}
Consequently, $\displaystyle n+\lceil \frac{m+1}2\rceil$ is the
smallest possible dimension of a compact absolute retract with the
$\DDP[m]{n}$-property.
\end{corollary}

Theorem~\ref{dim} implies that $\dim(X)\ge m+1$ for each metrizable
separable space $X$ with the $\HDDP[(2m+1)]{0}$-property. A similar
inequality also holds also for cohomological dimension.

\begin{theorem}\label{dim2} Let $X$ be a locally compact metrizable
$\LC[m]$-space having the $\HDDP[(2m+1)]{0}$-property. Then
$\dim_GX\ge m+1$ for any non-trivial Abelian group $G$.
\end{theorem}

In some cases the condition $X\in\HDDP[(2m+1)]{0}$ from Theorem
\ref{dim2} can be weakened to  $X\in\HDDP[2m]{0}$.

\begin{theorem}\label{dim2a}
Let $X$ be a locally compact $\LC[2m]$-space with the
$\HDDP[2m]{0}$-property and let $G$ be a non-trivial Abelian group.
The inequality $\dim_GX\ge m+1$ holds in each of the following
cases:
\begin{enumerate}
\item $G$ fails to be both divisible and periodic;
\item $G$ is a field;
\item $X$ is an $\ANR$-space.
\end{enumerate}
\end{theorem}

Theorem~\ref{dim2a} implies the following estimation for the
extension dimension of spaces $X\in\HDD[m]{0,0}$:

\begin{theorem}\label{extdim} Let  $X$ be a locally compact $\LC[m]$-space such that
$\edim X\le L$ for some $CW$-complex $L$. If $X\in\HDDP[m]{0}$, then
we have:
\begin{enumerate}
\item The homotopy groups $\pi_i(L)$ are trivial for all
$\displaystyle i<\frac{m}2$;
\item For $n=\lfloor {m}/2\rfloor$ the group
$\pi_n(L)$ is both divisible and periodic, and $\pi_n(L)=\widetilde
H_n(L)$;
\item $\pi_i(L)=0$ for all $\displaystyle i\le \frac{m}2$
provided $X$ is an ANR-space.
\end{enumerate}
\end{theorem}

Finally, we discuss the dimension properties of spaces $X\in\HDD[\infty]{0,0}$.

\begin{theorem} Let $X$ be a locally compact metrizable
$\LC[\infty]$-space with the $\HDD[\infty]{0,0}$-property. Then
\begin{enumerate}
\item All points of $X$ are homological $Z_\infty$-points;
\item $X$ fails to be $\trt$-dimensional;
\item If $\edim X\le L$ for some $CW$-complex $L$, then $L$ is contractible;
\item If $X$ is locally contractible, then $X$ is not a $C$-space.
\end{enumerate}
\end{theorem}

The first item of this theorem follows from Theorem~\ref{Zset}(7).
The last three items follow from Theorem~\ref{dimZn}(1) and
Theorem~\ref{edimZn}.

\subsection{Some Examples and Open Problems}
First, we discuss the problem of distinguishing between the
$\HDD[m]{n,k}$-properties for various $m,n,k$. Let us note that if
an Euclidean space $E$ has the $\HDD[m]{n,k}$-property for some
$m,n,k$, then $E$ has the $\HDD[a]{b,c}$-property for all
non-negative integers $a,b,c$ with $a+b+c\le n+m+k$. This feature is
specific for Euclidean spaces and does not hold in the general case.
For example, each dendrite $D$ with a dense set of end-points has
the $\HDD[0]{0,2}$-property (and in fact, $\HDD[0]{0,\infty}$) but
doesn't have the $\HDD[0]{1,1}$-property. Next example from
Daverman's book \cite{Dav} shows that the properties $\HDD[0]{0,2}$
and $\HDD[0]{1,1}$ are completely incomparable.

\begin{example}\label{ex1} There is a 2-dimensional absolute retract
$\Lambda\subset \IR^3$ with $\HDD[0]{1,1}$-property that fails to have the
$\HDD[0]{0,2}$-property.
\end{example}

\begin{question} Does the space $\Lambda$ from Example~\ref{ex1} possess the
$\HDD[2]{0,0}$-property$?$
\end{question}

It follows from Theorem~\ref{division}(1) that a Polish
$\LC[\infty]$-space $X$ has the $\HDD[\infty]{0,0}$-property
provided $X\times\IR\in\HDD[\infty]{0,0}$ and $X\in\HDD[2]{0,0}$. We
do not know if the latter condition is essential.

\begin{question} Does a compact absolute retract $X$ possess the
$\HDD[\infty]{0,0}$-property provided $X\times\mathbb I$ has that
property$?$ (Let us observe that $X\times\mathbb
I\in\HDD[\infty]{0,0}$ implies $X\times\mathbb
I\in\HDD[\infty]{1,\infty}$).
\end{question}

This question is equivalent to another intriguing one

\begin{question} Does a compact absolute retract $X$ contain a
$Z_2$-point provided all points of the product $X\times\mathbb I$
are $Z_\infty$-points$?$
\end{question}

\begin{problem} Let $X$ be a compact $\AR$ with the $\HDD[\infty]{0,0}$-property.
\begin{enumerate}
\item Is there any $Z_2$-point in $X?$
\item Is $X$ strongly infinite-dimensional $?$
\item Is $X\times\mathbb I$ homeomorphic to the Hilbert cube $?$
\end{enumerate}
\end{problem}

\begin{problem} Is a space $X\in\HDD[0]{2,2}$ homeomorphic to the
Hilbert cube $Q$ provided some finite power of $X$ is homeomorphic
to $Q$\textup{?}
\end{problem}

There are three interesting examples relevant to these questions.
The first of them was constructed by Singh in \cite{sin}, the second
by Daverman and Walsh in \cite{DW} and the third by Banakh and
Repov\v s in \cite{BR}.

\begin{example}[Singh]\label{ex2} There is a space $X$ possessing the following properties:
\begin{enumerate}
\item $X$ is a compact absolute retract;
\item $X$ contains no topological copy of the 2-disk $\mathbb I^2$;
\item $X\times\mathbb I$ is homeomorphic to the Hilbert cube;
\item All points of $X$ except for countably many are $Z_2$-points;
\item $X\in\HDD[\infty]{0,0}$;
\item $X\notin \HDD[2]{0,2}\cup \HDD[0]{2,2}$;
\item $X\times\mathbb I\in \HDD[\infty]{\infty,\infty}$.
\end{enumerate}
\end{example}

\begin{example}[Daverman-Walsh]\label{ex3} There is a space $X$ possessing the following properties:
\begin{enumerate}
\item $X$ is a compact absolute retract;
\item $X\times\mathbb I$ is homeomorphic to the Hilbert cube;
\item each point of $X$ is a $Z_\infty$-point;
\item $X\in\HDD[\infty]{0,\infty}\cap \HDD[0]{1,\infty}$
\item $X\notin \HDD[0]{2,2}$;
\item $X\times\mathbb I\in \HDD[\infty]{\infty,\infty}$.
\end{enumerate}
\end{example}

\begin{example}[Banakh-Repov\v s] There is a countable family
$\mathcal X$ of spaces such that
\begin{enumerate}
\item the product $X\times Y$ of any two different spaces
$X,Y\in\mathcal X$ is homeomorphic to the Hilbert cube;
\item no
finite power $X^k$ of any space $X\in\mathcal X$  is homeomorphic to
$Q$.
\end{enumerate}
\end{example}

It is interesting to note that there is no uncountable family
$\mathcal X$ possessing the properties (1) and (2) from Example 4,
see \cite{BR}.

It may be convenient to describe the $\HDD[m]{n,k}$-properties of a
space $X$ using the following sets $$
\begin{aligned}
\HDD[*]{*,*}(X)&=\{(m,n,k)\in\w^3:\mbox{$X$ has the
$\HDD[m]{n,k}$-property}\},\\
\HDD[0]{*,*}(X)&=\{(n,k)\in\w^2:\mbox{$X$ has the
$\HDD[0]{n,k}$-property}\},\\
\HDDP[*]{*}(X)&=\{(m,n)\in\w^2:\mbox{$X$ has the
$\HDDP[m]{n}$-property}\}.
\end{aligned}
$$

\begin{problem} Describe the geometry of the sets $\HDD[*]{*,*}(X)$,
$\HDD[0]{*,*}(X)$ and $\HDDP[*]{*}(X)$ for a given space $X$. Which
subsets of $\w^3$ or $\w^2$ can be realized as the sets
$\HDD[*]{*,*}(X)$, $\HDD[0]{*,*}(X)$ or $\HDDP[*]{*}(X)$ for a
suitable $X$\textup{?}
\end{problem}

In fact, we can consider the following partial pre-order
$\HDDmore$ on the set $\w^3$:\\
\mbox{$(m,n,k)\HDDmore(a,b,c)$} if each space $X$ with the
$\HDD[m]{n,k}$-property has also the $\HDD[a]{b,c}$-property.

\begin{problem} Describe the properties of the partial preorder
$\HDDmore$ on $\w^3$.
\end{problem}

By Proposition~\ref{lefschetz}(2), a paracompact space $X$ is
Lefschetz $\ANE[n]$ for a finite $n$ if and only if $X$ is an
$\LC[n-1]$-space. Consequently, the product of two paracompact
$\ANE[n]$-spaces is an $\ANE[n]$-space for every finite $n$.

\begin{problem} Is the product of two $($paracompact$)$ Lefschetz $\ANE[\infty]$-spaces a Lefschetz $\ANE[\infty]$-space?
\end{problem}



\end{document}